\documentclass[a4, 10pt]{amsart}
\usepackage{amssymb}
\usepackage{amstext}
\usepackage{amsmath}
\usepackage{amscd}
\usepackage{latexsym}
\usepackage{amsfonts}

\theoremstyle{plain}
\newtheorem{thm}{Theorem}[section]
\newtheorem*{thm*}{Theorem}
\newtheorem*{cor*}{Corollary}

\newtheorem{prop}[thm]{Proposition}
\newtheorem{lem}[thm]{Lemma}

\newtheorem{cor}[thm]{Corollary}

\newtheorem*{claim*}{Claim}

\theoremstyle{definition}
\newtheorem{defn}[thm]{Definition}
\newtheorem{ex}[thm]{Example}
\newtheorem{rem}[thm]{Remark}
\newtheorem{conj}[thm]{Conjecture}
\newtheorem{ques}[thm]{Question}

\theoremstyle{remark}
\newtheorem*{pf}{{\sl Proof}}

\newtheorem*{cpf}{{\sl Proof of Claim}}
\newtheorem*{lrpf}{{\sl Proof of Proposition \ref{laprap}}}

\numberwithin{equation}{thm}
\def\Hom{\mathrm{Hom}}
\def\lhom{\underline{\Hom}}

\def\Ext{\mathrm{Ext}}
\def\cext{\mathrm{\widehat{Ext}}\mathrm{}}

\def\Mod{\mathrm{Mod}}
\def\mod{\mathrm{mod}}
\def\Coker{\mathrm{Coker}}

\def\Im{\mathrm{Im}}
\def\tr{\mathrm{Tr}}

\def\m{\mathfrak m}

\def\p{\mathfrak p}
\def\q{\mathfrak q}
\def\P{\mathfrak P}

\def\Z{\Bbb Z}

\def\H{\mathrm{H}}

\def\Gdim{\mathrm{G}\mathrm{dim}}

\def\Kdim{\mathrm{dim}}
\def\depth{\mathrm{depth}}
\def\Supp{\mathrm{Supp}}

\def\Ass{\mathrm{Ass}}

\def\Min{\mathrm{Min}}

\def\id{\mathrm{id}}

\def\grade{\mathrm{grade}}
\def\codim{\mathrm{codim}}
\def\Spec{\mathrm{Spec}}

\def\lap{\mathrm{lap}}
\def\rap{\mathrm{rap}}
\def\A{{\mathcal A}}

\def\C{{\mathcal C}}
\def\F{{\mathcal F}}
\def\X{{\mathcal X}}
\def\Y{{\mathcal Y}}
\def\W{{\mathcal W}}
\def\G{{\mathcal G}}
\def\M{{\mathcal M}}
\def\lg{\underline{\G}}
\def\gt{\G ^{\bot}}
\def\tgt{{}^{\bot} (\G ^{\bot})}

\begin{document}

\title{Remarks on modules approximated by G-projective modules}
\author{Ryo Takahashi}
\address{Department of Mathematics, School of Science and Technology, Meiji University, Kawasaki 214-8571, Japan}
\email{takahasi@math.meiji.ac.jp}
\thanks{
\endgraf
}
\maketitle
\begin{abstract}
Let $R$ be a commutative Noetherian Henselian local ring.
Denote by $\mod\,R$ the category of finitely generated $R$-modules, and by $\G$ the full subcategory of $\mod\,R$ consisting of all G-projective $R$-modules.
In this paper, we consider when a given $R$-module has a right $\G$-approximation.
For this, we study the full subcategory $\rap\,\G$ of $\mod\,R$ consisting of all $R$-modules that admit right $\G$-approximations.
We investigate the structure of $\rap\,\G$ by observing $\G$, $\gt$ and $\lap\,\G$, where $\lap\,\G$ denotes the full subcategory of $\mod\,R$ consisting of all $R$-modules that admit left $\G$-approximations.
On the other hand, we also characterize $\rap\,\G$ in terms of Tate cohomologies.
We give several sufficient conditions for $\G$ to be contravariantly finite in $\mod\,R$.
\end{abstract}
\tableofcontents
\section{Introduction}

In the 1960's, Auslander \cite{Auslander} defined a homological invariant for finitely generated modules which he called Gorenstein dimension, G-dimension for short.
The value of G-dimension ranges from zero to infinity, and modules of finite G-dimension enjoy a lot of nice properties; they behave similarly to finitely generated modules over Gorenstein local rings.
Moreover, modules of finite G-dimension are resolved into finitely generated G-projective modules.
Thus, the class of finitely generated G-projective modules plays an essential role in considering G-dimension.
In this paper, we will observe finitely generated G-projective modules, and study the behavior of the class of them, which will be denoted by $\G$, in the category of finitely generated modules.
The main purpose of this paper is to know when a given module is approximated by the finitely generated G-projective modules.

Throughout the present paper, $R$ denotes a commutative Noetherian Henselian local ring with maximal ideal $\m$ and residue class field $k$, and all $R$-modules are assumed to be finitely generated modules.
We denote by $\mod\,R$ the category of finitely generated $R$-modules.
By a {\it subcategory} of $\mod\,R$ we always mean a full subcategory which is closed under isomorphisms.
(Recall that a subcategory $\X$ of $\mod\,R$ is said to be closed under isomorphisms provided that for any two objects $M, N$ of $\mod\,R$, if $M$ belongs to $\X$ and $N$ is isomorphic to $M$ then $N$ also belongs to $\X$.)
Similarly, a {\it subcategory} of a subcategory $\X$ of $\mod\,R$ always means a full subcategory of $\X$ which is closed under isomorphisms.

It is a well-known result due to Auslander and Buchweitz \cite{ABu} that if $R$ is Cohen-Macaulay, then for each $R$-module $M$, there exists a short exact sequence
$$
0 \to Y \to X \overset{f}{\to} M \to 0
$$
of $R$-modules such that $X$ is maximal Cohen-Macaulay and $Y$ is of finite injective dimension.
Such an exact sequence is called a {\it Cohen-Macaulay approximation} of $M$.
The reason why this is called an approximation is based on the fact that any homomorphism from a maximal Cohen-Macaulay $R$-module to $M$ factors through the homomorphism $f$ in the exact sequence.
In general, for a subcategory $\X$ of $\mod\,R$, a homomorphism $f:X\to M$ of $R$-modules with $X\in\X$ is called a right $\X$-approximation of $M$ if any homomorphism $f':X'\to M$ with $X'\in\X$ factors through $f$.
If any $R$-module in $\mod\,R$ has a right $\X$-approximation, then $\X$ is said to be contravariantly finite in $\mod\,R$.

Let $\G$ denote the subcategory of $\mod\,R$ which consists of all G-projective $R$-modules.
It is known that over a Gorenstein local ring, a finitely generated module is G-projective if and only if it is maximal Cohen-Macaulay.
Hence it follows from the above result of Auslander and Buchweitz that if $R$ is Gorenstein, then $\G$ is contravariantly finite in $\mod\,R$.
The author \cite{Takahashi3} conjectures that the converse also holds under a due assumption:

\begin{conj}
Suppose that there is a nonfree G-projective $R$-module.
If $\G$ is contravariantly finite in $\mod\,R$, then $R$ is Gorenstein.
\end{conj}

If this conjecture is true, then it holds that there exist infinitely many isomorphism classes of indecomposable G-projective $R$-modules whenever $R$ is non-Gorenstein and possesses a nonfree G-projective module.
Yoshino \cite[Theorem 6.1]{Yoshino3} proved that this conjecture is true for a certain Artinian local ring, and the author proved that it is true for any Henselian local ring of depth at most two; see Lemma \ref{mymain} below.
However, it is unknown whether the conjecture is true for a local ring of depth more than two.

In the present paper, in connection with this problem, we will consider $R$-modules having right $\G$-approximations; we want to give as many conditions as possible for a given $R$-module to have a right $\G$-approximation.
For this, we will observe such $R$-modules from various points of view.
Several subcategories of $\mod\,R$ which are associated to $\G$ will be introduced and studied.

Firstly, the subcategories $\gt$ and $\tgt$ of $\mod\,R$ will appear.
The former consists of all $R$-modules $Y$ such that $\Ext _R ^i (X,Y)=0$ for all $X\in\G$ and $i>0$, and the latter consists of all $R$-modules $Z$ such that $\Ext _R ^i (Z, Y)=0$ for all $Y\in\gt$ and $i>0$.
The subcategory $\gt$ is thick, namely, for an exact sequence $0 \to L \to M \to N \to 0$ of $R$-modules, if two of $L,M,N$ belong to $\gt$, then so does the third.
The subcategory $\tgt$ contains $\G$, and we will prove that $\tgt$ coincides with $\G$ if $R$ is a generically Gorenstein Cohen-Macaulay local ring admitting a canonical module.
After that, over such a ring, it will be shown that $\G$ is contravariantly finite in $\C$ (hence in $\mod\,R$) if $\gt\cap\C$ is covariantly finite in $\C$, where $\C$ denotes the subcategory of $\mod\,R$ consisting of all maximal Cohen-Macaulay $R$-modules.

Secondly, the subcategory $\rap\,\G$ of $\mod\,R$ will appear.
This subcategory consists of all $R$-modules that have right $\G$-approximations.
We shall prove that $\rap\,\G$ is a thick subcategory of $\mod\,R$, and is the smallest subcategory of $\mod\,R$ containing $\G$ and $\gt$ which is closed under direct summands and extensions.
As a corollary, one can prove that if $R$ is a non-Gorenstein local ring of depth at most two and there is a nonfree G-projective $R$-module, then no syzygy of the $R$-module $k$ admits a right $\G$-approximation.
Moreover, the fact that any module of finite G-dimension admits a right $\G$-approximation is obtained immediately.
We shall also show that $\G$ is contravariantly finite in $\mod\,R$ if $R$ is reduced and $\rap\,\G$ contains $\lap\,G$, which denotes the subcategory of $\mod\,R$ consisting of all $R$-modules having left $\G$-approximations.

On the other hand, we will give a criterion for a given $R$-module to have a right $\G$-approximation, in terms of Tate cohomologies.
To be concrete, we shall prove that the condition that an $R$-module $M$ has a right $\G$-approximation is equivalent to finite generation, finite presentation, and projectivity of $\cext _R ^i (-,M)|_{\lg}$ in the functor category of $\lg$ for some/any integer $i$, where $\lg$ denotes the stable category of $\G$.

In this paper, we will often refer to the papers \cite{AR} and \cite{AS}, which deal with modules over artin algebras.
Since the proofs of the results in those papers (to which we will refer) are completely categorical in nature, they carry over verbatim to the context of Henselian local rings.

We end this section by recalling the definitions of several conditions on a subcategory of $\mod\,R$ which we will often use in this paper.
For the definition of Auslander transpose, see the following part of Proposition \ref{gcm}.

\begin{defn}
For a subcategory $\X$ of $\mod\,R$, we say that
\begin{enumerate}
\item[{\rm (1)}]
$\X$ is {\it closed under finite (direct) sums} (resp. {\it closed under (direct) summands}) provided that for $M, N\in\mod\,R$ if $M, N\in\X$ then $M\oplus N\in\X$ (resp. if $M\oplus N\in\X$ then $M, N\in\X$).
\item[{\rm (2)}]
$\X$ is {\it closed under extensions} (resp. {\it closed under kernels of epimorphisms}, {\it closed under cokernels of monomorphisms}) provided that for any short exact sequence $0 \to L \to M \to N \to 0$ in $\mod\,R$, if $L, N\in\X$ then $M\in\X$ (resp. if $M, N\in\X$ then $L\in\X$, if $L, M\in\X$ then $N\in\X$).
\item[{\rm (3)}]
$\X$ is {\it closed under syzygies} (resp. {\it closed under (Auslander) transposes}) for any $X\in\X$ one has $\Omega X\in\X$ (resp. {\it $\tr X\in\X$}).
\end{enumerate}
\end{defn}

\section{$\F$-approximations}

In this section, we will mainly study the properties of right and left approximations of modules by free modules, which will be used in the later sections.
Before stating the definitions of right and left approximations, we recall the notions of right and left minimal homomorphisms which are introduced in \cite{AS}.

Let $\rho : M\to N$ be a homomorphism of $R$-modules.
We say that $\rho$ is {\it right minimal} if any endomorphism $f: M \to M$ satisfying $\rho = \rho f$ is an automorphism.
Dually, we say that $\rho$ is {\it left minimal} if any endomorphism $g: N\to N$ satisfying $\rho = g \rho$ is an automorphism.

\begin{defn}
Let $\X$ be a subcategory of $\mod\,R$.
\begin{enumerate}
\item[{\rm (1)}]
Let $\phi : X\to M$ be a homomorphism from $X\in\X$ to $M\in\mod\,R$.
\begin{enumerate}
\item[{\rm (i)}]
We call $\phi$ or $X$ a {\it right $\X$-approximation} of $M$ if for any homomorphism $\phi ': X'\to M$ with $X'\in\X$ there exists a homomorphism $f: X'\to X$ such that $\phi '=\phi f$.
\item[{\rm (ii)}]
Assume that $\phi$ is a right $\X$-approximation of $M$.
We call $\phi$ or $X$ a {\it minimal right $\X$-approximation} of $M$ if $\phi$ is right minimal.
\end{enumerate}
\item[{\rm (2)}]
Let $\phi : M\to X$ be a homomorphism from $M\in\mod\,R$ to $X\in\X$.
\begin{enumerate}
\item[{\rm (i)}]
We call $\phi$ or $X$ a {\it left $\X$-approximation} of $M$ if for any homomorphism $\phi ': M\to X'$ with $X'\in\X$ there exists a homomorphism $f: X\to X'$ such that $\phi '=f\phi$.
\item[{\rm (ii)}]
Assume that $\phi$ is a left $\X$-approximation of $M$.
We call $\phi$ or $X$ a {\it minimal left $\X$-approximation} of $M$ if $\phi$ is left minimal.
\end{enumerate}
\end{enumerate}
\end{defn}

A right $\X$-approximation (resp. minimal right $\X$-approximation, left $\X$-approximation, minimal left $\X$-approximation) is also called a {\it $\X$-precover} (resp. {\it $\X$-cover}, {\it $\X$-preenvelope}, {\it $\X$-envelope}).
It is easily seen by definition that a minimal right (resp. left) $\X$-approximation is uniquely determined up to isomorphism whenever it exists.
For a subcategory $\X$ of $\mod\,R$ closed under direct summands, an $R$-module having a right (resp. left) $\X$-approximation also has a minimal right (resp. left) $\X$-approximation; see \cite[Proposition 2.4]{Takahashi3}.

\begin{defn}
Let $\M$ be a subcategory of $\mod\,R$, and let $\X$ be a subcategory of $\M$.
Then we say that $\X$ is {\it contravariantly finite} (resp. {\it covariantly finite}) in $\M$ if any $R$-module in $\M$ has a right (resp. left) $\X$-approximation.
If $\X$ is both covariantly finite and contravariantly finite in $\M$, then $\X$ is said to be {\it functorially finite} in $\M$.
\end{defn}

A contravariantly finite (resp. covariantly finite) subcategory is also called a {\it precovering} (resp. {\it preenveloping}) subcategory.

We denote by $\F$ the subcategory of $\mod\,R$ consisting of all free $R$-modules.
From now on, we shall consider right and left $\F$-approximations.
Recall that a homomorphism $f : M \to N$ of $R$-modules is said to be {\it minimal} if the induced homomorphism $f\otimes _R k : M\otimes _R k \to N\otimes _R k$ is an isomorphism.
Note from Nakayama's lemma that every minimal homomorphism is surjective.
Let $\nu _R (M)$ denote the minimal number of generators of an $R$-module $M$, i.e., $\nu _R (M)=\Kdim _k (M\otimes _R k)$.
Set $(-)^{\ast}=\Hom _R (-, R)$.
The following result is easily obtained.

\begin{prop}\label{freecover}
Let $M$ be an $R$-module.
\begin{enumerate}
\item[{\rm (1)}]
Let $\phi : R^n \to M$ be a homomorphism of $R$-modules.
The following conditions are equivalent:
\begin{enumerate}
\item[{\rm (i)}]
$\phi$ is a minimal right $\F$-approximation of $M$;
\item[{\rm (ii)}]
$\phi$ is surjective and $n=\nu _R (M)$.
\end{enumerate}
\item[{\rm (2)}]
Let $f_1, f_2, \dots, f_n$ be a minimal system of generators of $M^{\ast}$.
Then the homomorphism
$$
f=\left(
\begin{smallmatrix}
f_1\\
\vdots\\
\\
f_n
\end{smallmatrix}
\right) : M \to R^n
$$
is a minimal left $\F$-approximation of $M$.
\item[{\rm (3)}]
Let $\sigma : M \to M^{\ast\ast}$ be the natural homomorphism and $\phi : F \to M^{\ast}$ a minimal right $\F$-approximation.
Then the composite map $\phi ^{\ast} \sigma: M \to F^{\ast}$ is a minimal left $\F$-approximation.
\end{enumerate}
\end{prop}

An $R$-module $M$ is said to be {\it torsionless} (resp. {\it reflexive}) if the natural homomorphism $M \to M^{\ast\ast}$ is injective (resp. bijective).
Here we state a property of left $\F$-approximations of torsionless modules.

\begin{prop}\label{envelope3}
The following are equivalent:
\begin{enumerate}
\item[{\rm (1)}]
$M$ is torsionless;
\item[{\rm (2)}]
Every left $\F$-approximation of $M$ is an injective homomorphism;
\item[{\rm (3)}]
Some left $\F$-approximation of $M$ is an injective homomorphism.
\end{enumerate}
\end{prop}

\begin{pf}
Note by \cite[Lemma 3.4]{EG} that an $R$-module is torsionless if and only if it is a first syzygy.
Let $\psi : M \to R^n$ be a left $\F$-approximation of $M$.
If $M$ is torsionless, then there is an injective homomorphism $\rho : M \to R^m$.
The definition of a left approximation says that $\rho$ factors through $\psi$, which shows that $\psi$ is also an injective homomorphism.
\qed
\end{pf}

Note from the above proposition that a minimal left $\F$-approximation is not necessarily an injective homomorphism.

Let $M$ be an $R$-module.
Take its minimal right $\F$-approximation $\pi : F \to M$.
The {\it first syzygy} $\Omega M=\Omega ^1 M$ of $M$ is defined as the kernel of the homomorphism $\pi$ (cf. Proposition \ref{freecover}(1)), and the {\it $n$th syzygy} $\Omega ^n M$ of $M$ is defined inductively: $\Omega ^n M=\Omega (\Omega ^{n-1} M)$ for $n\geq 2$.
Dually to this, we can define the cosyzygies of a given $R$-module.

\begin{defn}
Let $M$ be an $R$-module.
\begin{enumerate}
\item[$(1)$]\ Take the minimal left $\F$-approximation $\theta : M \to F$ of $M$.
We set $\Omega ^{-1} M = \Coker\,\theta$, and call it the {\it first cosyzygy} of $M$.
\item[$(2)$]\ Let $n\geq 2$.
Assume that the $(n-1)$th cosyzygy $\Omega ^{-(n-1)} M$ is defined.
Then we set $\Omega ^{-n} M = \Omega ^{-1} (\Omega ^{-(n-1)} M)$ and call it the {\it $n$th cosyzygy} of $M$.
\end{enumerate}
\end{defn}

An $R$-module is said to be {\it stable} if it has no nonzero free $R$-summand.
Cosyzygies are always stable:

\begin{prop}\label{stable}
For any $R$-module $M$ and any positive integer $n$, the $R$-module $\Omega ^{-n} M$ is stable.
\end{prop}

\begin{pf}
We have only to show that $\Omega ^{-1} M$ is stable.
Denote by $\theta : M\to R^m$ the minimal left $\F$-approximation of $M$.
There is an exact sequence
$$
M \overset{\theta}{\to} R^m \to \Omega ^{-1} M \to 0.
$$
Suppose that $\Omega ^{-1} M$ is not stable.
Then there exists a surjective homomorphism $\varepsilon : \Omega ^{-1} M \to R$.
We can write a commutative diagram
$$
\begin{CD}
@. M @>{\theta}>> R^m @>>> \Omega ^{-1} M @>>> 0\\
@. @V{\theta '}VV @| @V{\varepsilon}VV\\
0 @>>> R^{m-1} @>{f}>> R^m @>>> R @>>> 0
\end{CD}
$$
with exact rows.
Since $f$ is a split monomorphism, there is a homomorphism $g : R^m \to R^{m-1}$ such that $gf=1$.
Noting that $\theta = f\theta '$, we have $fg\theta =\theta$.
Hence $fg$ is an automorphism because $\theta$ is a minimal left $\F$-approximation.
Thus the homomorphism $f: R^{m-1} \to R^m$ must be surjective.
But this is a contradiction, which proves that $\Omega ^{-1} M$ is stable.
\qed
\end{pf}

For a subcategory $\X$ of $\mod\,R$, we denote by $\X ^{\sf L}$ (resp. ${}^{\sf L}\X$) the subcategory of $\mod\,R$ consisting of all $R$-modules $M$ such that $\Ext _R ^1 (X, M)=0$ (resp. $\Ext _R ^1 (M, X)=0$) for all $X\in\X$.
The proposition below follows from a Wakamatsu's lemma \cite[Lemma 2.1.2]{Xu}.

\begin{prop}\label{cosyz}
Any cosyzygy belongs to ${}^{\sf L}\F$, namely
$$
\Ext _R ^1 (\Omega ^{-1}M, R)=0
$$
for any $R$-module $M$.
\end{prop}

\section{Basic properties of $\G$}

In this section, we will study several basic properties of a G-projective module and G-dimension.
Let us recall their definitions.

\begin{defn}
Denote by $(-)^{\ast}$ the $R$-dual functor $\Hom _R (-, R)$.
\begin{enumerate}
\item[{\rm (1)}]
We say that an $R$-module $X$ is {\it G-projective} if the following three conditions hold:
\begin{enumerate}
\item[{\rm (i)}]
The natural homomorphism $X\to X^{\ast\ast}$ is an isomorphism,
\item[{\rm (ii)}]
$\Ext _R ^i (X, R)=0$ for any $i>0$,
\item[{\rm (iii)}]
$\Ext _R ^i (X^{\ast}, R)=0$ for any $i>0$.
\end{enumerate}
We denote by $\G$ the full subcategory of $\mod\,R$ consisting of all G-projective $R$-modules.
\item[{\rm (2)}]
Let $M$ be an $R$-module.
If there exists an exact sequence
$$
0 \to X_n \to X_{n-1} \to \cdots \to X_1 \to X_0 \to M \to 0
$$
of $R$-modules with $X_i\in\G$ for each $i$, then we say that $M$ has {\it G-dimension at most $n$}, and write $\Gdim _R\, M\leq n$.
If such an integer $n$ does not exist, then we say that $M$ has {\it infinite G-dimension}, and write $\Gdim _R\, M=\infty$.
\end{enumerate}
\end{defn}

If an $R$-module $M$ has G-dimension at most $n$ but does not have G-dimension at most $n-1$, then we say that $M$ has {\it G-dimension $n$}, and write $\Gdim _R\, M=n$.
Note that being G-dimension zero is equivalent to being G-projective.

The result below immediately follows from the Auslander-Bridger formula \cite[(1.4.8)]{Christensen}.

\begin{prop}\label{gcm}
Let $R$ be a Cohen-Macaulay local ring.
Then every G-projective $R$-module is maximal Cohen-Macaulay.
\end{prop}

Let
$$
F_1 \overset{\delta}{\to} F_0 \to M \to 0
$$
be a minimal free presentation of an $R$-module $M$, that is to say, it is an exact sequence such that $F_0, F_1$ are free $R$-modules and the image of the homomorphism $\delta$ is contained in $\m F_0$.
We denote by $\tr M$ the cokernel of the $R$-dual homomorphism ${\delta}^{\ast} : F_0^{\ast} \to F_1^{\ast}$.
It is called the {\it (Auslander) transpose} or {\it Auslander dual} of $M$.
We should note that the module $\tr M$ is uniquely determined up to isomorphism because we defined it by using a minimal free presentation of $M$.
We should also note that $M$ is isomorphic to $\tr (\tr M)$ up to free summand.
For more details on transposes, refer to \cite{ABr} and \cite{Masek}.

An $R$-complex
$$
F_{\bullet} = (\cdots \overset{d_2}{\longrightarrow} F_1 \overset{d_1}{\longrightarrow} F_0 \overset{d_0}{\longrightarrow} F_{-1} \overset{d_{-1}}{\longrightarrow} F_{-2} \overset{d_{-2}}{\longrightarrow} \cdots)
$$
is said to be a {\it complete resolution} of an $R$-module $M$ if the following three conditions hold:
\begin{enumerate}
\item[{\rm (a)}]
$F_i\in\F$ for any $i\in\Z$,
\item[{\rm (b)}]
$\H _i (F_{\bullet})=0=\H ^i ((F_{\bullet})^{\ast})$ for any $i\in\Z$,
\item[{\rm (c)}]
$\Im\,d_0=M$.
\end{enumerate}

We present properties of a G-projective module which we will often use later.

\begin{prop}\label{g}
\begin{enumerate}
\item[{\rm (1)}]
The following are equivalent for an $R$-module $M$:
\begin{enumerate}
\item[{\rm (i)}]
$M$ is G-projective;
\item[{\rm (ii)}]
$\Ext _R ^i (M, R) = 0 = \Ext _R ^i (\tr M, R)$ for every $i>0$;
\item[{\rm (iii)}]
$M$ admits a complete resolution.
\end{enumerate}
\item[{\rm (2)}]
If an $R$-module $M$ is G-projective, then so are $M^{\ast}$, $\tr M$, $\Omega M$ and $\Omega ^{-1}M$.
\end{enumerate}
\end{prop}

\begin{pf}
(1) This statement is proved in \cite[Proposition (3.8)]{ABr} and \cite[(4.1.4)]{Christensen}.

(2) Let $M$ be a G-projective $R$-module.
By definition, $M^{\ast}$ is G-projective.
The statement (1) shows that $\tr M$ is G-projective.
It follows from \cite[Lemma 2.3]{AM} that $\Omega M$ is G-projective.
Noting that $\Omega ^{-1}M$ is isomorphic to $(\Omega (M^{\ast}))^{\ast}$ by Proposition \ref{freecover}(3), one sees that $\Omega ^{-1}M$ is G-projective.
\qed
\end{pf}

Here we give a remark on the structure of a complete resolution; it consists of right and left $\F$-approximations:

\begin{prop}\label{compresol}
Let
$$
\cdots \overset{d_2}{\longrightarrow} F_1 \overset{d_1}{\longrightarrow} F_0 \overset{d_0}{\longrightarrow} F_{-1} \overset{d_{-1}}{\longrightarrow} F_{-2} \overset{d_{-2}}{\longrightarrow} \cdots
$$
be a complete resolution of an $R$-module $M$.
Let $\alpha : F_0 \to M$ be the surjective homomorphism induced by $d_0$, and let $\beta : M \to F_{-1}$ be the inclusion map.
Then $\alpha$ (resp. $\beta$) is a right (resp. left) $\F$-approximation of $M$.
\end{prop}

\begin{pf}
It is easy to see from the surjectivity of $\alpha$ that $\alpha$ is a right $\F$-approximation.
Take a free $R$-module $P$.
Noting by the definition of a complete resolution that $\Hom _R (F_{\bullet}, P)$ is an exact complex, one sees that the homomorphism
$$
\Hom _R (\beta , P): \Hom _R (F_{-1}, P) \to \Hom _R (M, P)
$$
is surjective, which means that $\beta$ is a left $\F$-approximation of $M$.
\qed
\end{pf}

\begin{defn}
A subcategory $\X$ of $\mod\,R$ is said to be {\it resolving} if the following hold:
\begin{enumerate}
\item[{\rm (1)}]
$\X$ contains $R$,
\item[{\rm (2)}]
$\X$ is closed under direct summands,
\item[{\rm (3)}]
$\X$ is closed under extensions,
\item[{\rm (4)}]
$\X$ is closed under kernels of epimorphisms.
\end{enumerate}
\end{defn}

For a given subcategory $\X$ of $\mod\,R$, we denote by $\X ^{\bot}$ (resp. ${}^{\bot}\X$) the subcategory of $\mod\,R$ consisting of all $R$-modules $M$ such that $\Ext _R ^i (X, M)=0$ (resp. $\Ext _R ^i (M, X)=0$) for all $X\in\X$ and $i>0$.
Also, we denote by $\widehat{\X}$ the subcategory of $\mod\,R$ consisting of all $R$-modules $M$ that have exact sequences
$$
0 \to X_n \to X_{n-1} \to \cdots \to X_1 \to X_0 \to M \to 0
$$
with $X_i\in\X$ for $0\leq i\leq n$.
Let $\Y$ be a subcategory of $\X$.
We say that $\Y$ is {\it Extinjective} in $\X$ if $\Y$ is contained in $\X ^{\bot}$.
We say that $\Y$ is a {\it cogenerator} for $\X$ if for any $X\in\X$ there exists an exact sequence $0 \to X \to Y \to X' \to 0$ with $Y\in\Y$ and $X'\in\X$.

Now, we can state a well-known result due to Auslander and Buchweitz.
For the proof, see \cite[Theorem 1.1, Proposition 3.6]{ABr}.

\begin{lem}[Auslander-Buchweitz]\label{abr}
Let $\X$ be a resolving subcategory of $\mod\,R$ with Extinjective cogenerator $\W$.
Then the following hold.
\begin{enumerate}
\item[{\rm (1)}]
$\X$ is contravariantly finite in $\widehat{\X}$,
\item[{\rm (2)}]
$\widehat{\W}=\X^{\bot}\cap\widehat{\X}$.
\end{enumerate}
\end{lem}

The subcategory $\G$ of $\mod\,R$ satisfies the assumptions of the above result:

\begin{prop}\label{gra}
$\G$ is a resolving subcategory of $\mod\,R$ with Extinjective cogenerator $\F$.
\end{prop}

\begin{pf}
It follows from \cite[Lemma 2.3]{AM} that $\G$ is a resolving subcategory of $\mod\,R$, and it is obvious from definition that $\F$ is Extinjective in $\G$.
Hence we have only to show that $\F$ is a cogenerator for $\G$.
Let $X\in\G$.
Then, since $X$ is torsionless by definition, Proposition \ref{envelope3} implies that one has an exact sequence $0 \to X \to F \to \Omega ^{-1}X \to 0$ with $F\in\F$.
According to Proposition \ref{g}(2), the module $\Omega ^{-1}X$ belongs to $\G$.
Thus $\F$ is a cogenerator for $\G$.
\qed
\end{pf}

Lemma \ref{abr}, Proposition \ref{gra} and \cite[(1.4.9)]{Christensen} give the following connections between G-projective modules and modules of finite G-dimension (see also \cite[Theorem 8.6]{AM}).

\begin{cor}\label{abg}
\begin{enumerate}
\item[{\rm (1)}]
Any $R$-module of finite G-dimension has a right $\G$-approximation.
\item[{\rm (2)}]
An $R$-module $M$ belongs to $\G ^{\bot}$ and has finite G-dimension if and only if $M$ has finite projective dimension.
\item[{\rm (3)}]
If $R$ is Gorenstein, then $\G$ is contravariantly finite in $\mod\,R$.
\end{enumerate}
\end{cor}

\section{The relationship between $\G$ and $\tgt$}

In this section, we will study the inclusion relation between the subcategories $\G$ and $\tgt$ of $\mod\,R$.
To be concrete, we will give several sufficient conditions for the subcategory $\G$ to coincide with the subcategory $\tgt$.

\begin{prop}\label{resappr}
Let $\G$ be a resolving subcategory of $\mod\,R$.
Then the following hold.
\begin{enumerate}
\item[{\rm (1)}]
One has $\G ^{\sf L} = \G ^{\bot}$.
\item[{\rm (2)}]
Suppose that an $R$-module $M$ has a right $\G$-approximation.
Then there exists an exact sequence
$$
0 \to Y \to X \to M \to 0
$$
of $R$-modules with $X\in\G$ and $Y\in\G ^{\bot}$.
\end{enumerate}
\end{prop}

\begin{pf}
The assertions actually hold for an arbitrary resolving subcategory of $\mod\,R$; see \cite[Lemma 3.2(a), Proposition 3.3(c)]{AR}.
Proposition \ref{gra} says that $\G$ is a resolving subcategory of $\mod\,R$.
\qed
\end{pf}

Using Proposition \ref{resappr}(2) and Corollary \ref{abg}(3), we easily obtain the following result.
(The proof is similar to that of \cite[Proposition 3.3(b)]{AR}.)

\begin{cor}\label{tgtcontravfin}
If $\G$ is contravariantly finite in $\tgt$, then $\G =\tgt$.
In particular, one has $\G = \tgt$ whenever $R$ is Gorenstein.
\end{cor}

Similarly, one can show that if $\G$ is contravariantly finite in ${}^{\sf L} (\G ^{\sf L})$ then $\G = {}^{\sf L} (\G ^{\sf L})$.

Next, we want to investigate the subcategory $\gt$ of $\mod\,R$.
Before that, let us recall the definition of a thick subcategory.

\begin{defn}
Let $\X$ be a subcategory of $\mod\,R$ which is closed under direct summands.
We say that $\X$ is {\it thick} provided that for any exact sequence $0 \to L \to M \to N \to 0$ in $\mod\,R$, if two of $L, M, N$ belong to $\X$ then so does the third.
\end{defn}

\begin{prop}\label{gtthick}
$\gt$ is a thick subcategory of $\mod\,R$ containing $R$.
In particular, $\gt$ is a resolving subcategory of $\mod\,R$.
\end{prop}

\begin{pf}
It is immediately follows from definition that $\gt$ is closed under direct summands and contains $R$.
Let $0 \to L \to M \to N \to 0$ be an exact sequence in $\mod\,R$.
We easily observe that if $L, N\in\gt$ then $M\in\gt$, and that if $L, M\in\gt$ then $N\in\gt$.
Suppose that $M, N\in\gt$.
Then it is seen that $\Ext _R ^i (X, L)=0$ for any $X\in\G$ and $i\geq 2$.
Fix $X\in\G$.
By Proposition \ref{envelope3}, there is an exact sequence $0 \to X \to F \to \Omega ^{-1}X \to 0$ where $F$ is a free $R$-module, and it follows from this sequence that $\Ext _R ^1 (X, L)\cong\Ext _R ^2 (\Omega ^{-1}X, L)=0$ because $\Omega ^{-1}X\in\G$ by Proposition \ref{g}(2).
Consequently we have $\Ext _R ^i (X, L)=0$ for any $i\geq 1$ and $X\in\G$, that is to say, $L$ is in $\gt$.
Thus $\gt$ is a thick subcategory of $\mod\,R$.
\qed
\end{pf}

Let $n$ be a positive integer.
An $R$-module $M$ is called {\it $n$-torsionfree} if $\Ext _R ^i (\tr M, R)=0$ for any integer $i$ with $1\leq i\leq n$.
Now we can prove the main result of this section.

\begin{thm}\label{gtgt}
The following are equivalent:
\begin{enumerate}
\item[{\rm (1)}]
$\G = \tgt$;
\item[{\rm (2)}]
Every module in $\tgt$ is torsionless.
\end{enumerate}
\end{thm}

\begin{pf}
(1) $\Rightarrow$ (2): By definition, every G-projective module is reflexive, hence torsionless.

(2) $\Rightarrow$ (1): First of all, let us observe the following claim.

\begin{claim*}
If a stable $R$-module $M$ belongs to $\tgt$, then $M$ is isomorphic to $\Omega (\Omega ^{-1}M)$ and $\Omega ^{-1}M$ also belongs to $\tgt$.
\end{claim*}

\begin{cpf}
Since $M$ is torsionless by assumption, it is a first syzygy by \cite[Lemma 3.4]{EG}.
According to Proposition \ref{envelope3}, there is an exact sequence $0 \to M \to R^r \to \Omega ^{-1} M \to 0$.
From the exact sequence and the stability of $M$, we see that $M$ is isomorphic to $\Omega (\Omega ^{-1}M)$.
Also, for $Y\in\gt$ and $i\geq 2$, we have $\Ext _R ^i (\Omega ^{-1}M, Y)$ is isomorphic to $\Ext _R ^{i-1} (M, Y)$, and the latter Ext module is zero.
Hence
\begin{equation}\label{claimkey}
\Ext _R ^i (\Omega ^{-1}M, Y)=0\text{ for any }Y\in\gt\text{ and }i\geq 2.
\end{equation}

Fix $Y\in\gt$.
There is an exact sequence $0 \to \Omega Y \to R^s \to Y \to 0$, and we get an exact sequence
$$
\Ext _R ^1 (\Omega ^{-1}M, R^s) \to \Ext _R ^1 (\Omega ^{-1}M, Y) \to \Ext _R ^2 (\Omega ^{-1}M, \Omega Y).
$$
By virtue of Proposition \ref{cosyz}, we have $\Ext _R ^1 (\Omega ^{-1}M, R^s)=0$.
Proposition \ref{gtthick} implies that $\Omega Y\in\gt$, hence $\Ext _R ^2 (\Omega ^{-1}M, \Omega Y)=0$ by \eqref{claimkey}.
Therefore $\Ext _R ^1 (\Omega ^{-1}M, Y)=0$.
Thus, by \eqref{claimkey} again, $\Ext _R ^i (\Omega ^{-1}M, Y)=0$ for any $Y\in\gt$ and $i\geq 1$, which means that $\Omega ^{-1}M$ belongs to $\tgt$.
\qed
\end{cpf}

Now we shall prove that $\G$ coincides with $\tgt$.
It is clear that $\G$ is contained in $\tgt$.
Let us observe the converse inclusion relation.
Take $M\in\tgt$.
Note that the subcategory $\tgt$ is closed under direct summands and that $\G$ is closed under finite direct sums.
Hence, to show that $M$ belongs to $\G$, without loss of generality, we can assume that $M$ is a stable $R$-module.
Fix $n>0$.
From the above claim and Proposition \ref{stable}, one sees that $M\cong\Omega ^n (\Omega ^{-n}M)$ and $\Omega ^{-n}M\in\tgt$.
Since $\F\subseteq\gt$, we have $\tgt\subseteq {}^{\bot}\F$.
Hence $\Ext _R ^i (M, R)=0$ and $\Ext _R ^i (\Omega ^{-n}M, R)=0$ for any $i>0$.
It is seen by \cite[Proposition (2.26)]{ABr} that $\Omega ^i (\Omega ^{-n}M)$ is $i$-torsionfree for $1\leq i\leq n$.
Particularly, $M\cong\Omega ^n (\Omega ^{-n}M)$ is $n$-torsionfree.
Therefore $\Ext _R ^i (\tr M, R)=0$ for $1\leq i\leq n$.
Consequently, we have $\Ext _R ^i (\tr M, R)=0$ for any $i>0$, and thus the module $M$ belongs to $\G$ by Proposition \ref{g}(1).
\qed
\end{pf}

Recall that a local ring $R$ is said to be {\it generically Gorenstein} if $R_{\p}$ is a Gorenstein local ring for every $\p\in\Min\,R$.
Using the above theorem, one can find a relatively general class of local rings $R$ satisfying $\G=\tgt$:

\begin{cor}\label{gengorcm}
Let $R$ be a generically Gorenstein Cohen-Macaulay local ring with canonical module $\omega$.
Then $\G = \tgt$.
\end{cor}

\begin{pf}
Let $M$ be an $R$-module in $\tgt$.
Proposition \ref{gcm} says that every $R$-module in $\G$ is maximal Cohen-Macaulay, equivalently, the canonical module $\omega$ belongs to $\gt$.
Hence one has $\Ext _R ^i (M, \omega)=0$ for every $i>0$, equivalently, $M$ is a maximal Cohen-Macaulay $R$-module.
Therefore $\Ass\,M$ is contained in $\Ass\,R$.
It is seen by \cite[Lemma 3.4, Theorem 3.5]{EG} that $M$ is torsionless.
Thus the assertion follows from Theorem \ref{gtgt}.
\qed
\end{pf}

The results appearing in this section naturally lead us to a question:

\begin{ques}
Is it always true that $\G$ coincides with $\tgt$?
\end{ques}

We should note from Corollary \ref{tgtcontravfin} that if $\G$ does not coincide with $\tgt$ then $\G$ is not contravariantly finite in $\mod\,R$.

\section{Right $\G$-approximations over Cohen-Macaulay rings}

We denote by $\C$ the subcategory of $\mod\,R$ consisting of all maximal Cohen-Macaulay $R$-modules, i.e., $R$-modules $M$ satisfying $\depth _R\, M=\Kdim\, R$.
In this section, we will consider contravariant finiteness of $\G$ in $\C$ over a Cohen-Macaulay local ring $R$ admitting the canonical module.

First of all, we introduce the definitions of right and left functor categories.
Let $\A$ be an additive category.
The {\it right functor category} of $\A$, which is denoted by $\Mod\,\A$, is defined as the category having additive contravariant functors from $\A$ to the category of abelian groups as the objects, and natural transformations between such two functors as the morphisms.
An object of $\Mod\,\A$ is called a {\it right $\A$-module}.
For $F\in\Mod\,\A$, we say that $F$ is {\it finitely generated} if there exists an exact sequence
$$
\Hom _{\A} (-, X) \to F \to 0
$$
in $\Mod\,\A$.
We say that $F$ is {\it finitely presented} if there exists an exact sequence
$$
\Hom _{\A} (-, X_1) \to \Hom _{\A} (-, X_0) \to F \to 0
$$
in $\Mod\,\A$.
We denote by $\mod\,\A$ the full subcategory of $\Mod\,\A$ consisting of all finitely presented right $\A$-modules.

Also, the {\it left functor category} $\A\,\Mod$ of $\A$ is defined as the category of additive covariant functors from $\A$ to the category of abelian groups.
A {\it left $\A$-module}, a {\it finitely generated} left $\A$-module, a {\it finitely presented} left $\A$-module and the category $\A\,\mod$ are defined dually.

For a functor $F$ from $\mod\,R$ to itself and a subcategory $\X$ of $\mod\,R$, we denote by $F|_{\X}$ the restriction of $F$ to $\X$. 

\begin{lem}\label{tycont}
Let $\M$ be a resolving subcategory of $\mod\,R$, and let $\Y$ be a subcategory of $\M$ which is closed under extensions.
Suppose that $\Y$ is covariantly finite in $\M$.
Then ${}^{\sf L}\Y\cap\M$ is contravariantly finite in $\M$.
\end{lem}

\begin{pf}
According to \cite[Corollary 1.5]{AR}, it suffices to prove that the left $\Y$-module $\Ext _R ^1(M,-)|_{\Y}$ is finitely generated for any $M\in\M$.
Let $M\in\M$.
There is an exact sequence $0 \to \Omega M \to R^n \to M \to 0$.
From this exact sequence, we get a surjective morphism of functors
$$
\Hom _R (\Omega M, -)|_{\Y} \to \Ext _R ^1 (M, -)|_{\Y}.
$$
On the other hand, since $\Y$ is covariantly finite in $\M$ and $\Omega M\in\M$, there exists a left $\Y$-approximation $\Omega M \to Y$.
Noting the definition of a left approximation, we can make another surjective morphism of functors
$$
\Hom _R (Y, -)|_{\Y} \to \Hom _R (\Omega M, -)|_{\Y}.
$$
Splicing these two morphisms together, we get a surjective morphism of functors $\Hom _R (Y, -)|_{\Y} \to \Ext _R ^1 (M, -)|_{\Y}$, which says that the left $\Y$-module $\Ext _R ^1 (M, -)|_{\Y}$ is finitely generated, as desired.
\qed
\end{pf}

\begin{lem}\label{cateq}
Let $R$ be a Cohen-Macaulay local ring with canonical module $\omega$.
Then one has the following:
\begin{enumerate}
\item[{\rm (1)}]
$\tgt = {}^{\bot}(\gt\cap\C)$,
\item[{\rm (2)}]
${}^{\bot}(\gt\cap\C)\cap\C = {}^{\sf L}(\gt\cap\C)\cap\C$.
\end{enumerate}
\end{lem}

\begin{pf}
(1) Noting that $\gt\cap\C$ is contained in $\gt$, we see that $\tgt$ is contained in ${}^{\bot}(\gt\cap\C)$.
To observe the converse inclusion relation, take $M\in {}^{\bot}(\gt\cap\C)$.
Since $R$ belongs to $\gt\cap\C$, we have $\Ext _R ^i (M, R)=0$ for any $i>0$.
Fix $Y\in\gt$.
Then $\Omega ^n Y$ is a maximal Cohen-Macaulay $R$-module for $n\gg 0$.
Since $\gt$ is resolving by Proposition \ref{gtthick}, the module $\Omega ^n Y$ also belongs to $\gt$.
Hence $\Omega ^n Y$ belongs to $\gt\cap\C$, and therefore $\Ext _R ^i (M, \Omega ^n Y)=0$ for any $i>0$.
Thus we obtain isomorphisms
$$
\Ext _R ^i (M, Y) \cong \Ext _R ^{i+1} (M, \Omega Y) \cong \cdots \cong \Ext _R ^{i+n} (M, \Omega ^n Y)=0
$$
for any $i>0$, which says that $M$ belongs to $\tgt$.

(2) It is obvious that ${}^{\bot}(\gt\cap\C)\cap\C$ is contained in ${}^{\sf L}(\gt\cap\C)\cap\C$.
Let $M\in {}^{\sf L}(\gt\cap\C)\cap\C$ and $Y\in\gt\cap\C$.
We want to prove $\Ext _R ^i (M, Y)=0$ for every $i>0$.
Denote by $(-)^{\dag}$ the canonical dual functor $\Hom _R (-, \omega )$.
We have exact sequences
$$
0 \to \Omega ^{j+1}(Y^{\dag}) \to R^{n_j} \to \Omega ^j(Y^{\dag}) \to 0
$$
for $j\geq 0$.
Since $Y$ is maximal Cohen-Macaulay, so is $Y^{\dag}$, and so is $\Omega ^j (Y^{\dag})$ for any $j\geq 0$.
Applying $(-)^{\dag}$ to the above exact sequences, we get exact sequences
$$
0 \to Y_j \to \omega ^{n_j} \to Y_{j+1} \to 0,
$$
where $Y_j = (\Omega ^j (Y^{\dag}))^{\dag}$.
Noting that $\omega$ belongs to $\gt$ because $\G$ is contained in $\C$ by Proposition \ref{gcm}, we see that if $Y_j\in\gt$ then $Y_{j+1}\in\gt$.
Since $Y_0 \cong Y\in\gt$, an inductive argument shows that $Y_j\in\gt$, hence $Y_j\in\gt\cap\C$, for $j\geq 0$.
Therefore we obtain $\Ext _R ^1 (M, Y_j)=0$ for every $j\geq 0$.
Noting that $\Ext _R ^i (M, \omega)=0$ for $i>0$ because $M\in\C$, we have isomorphisms
$$
\Ext _R ^i (M, Y) \cong \Ext _R ^{i-1} (M, Y_1) \cong \cdots \cong \Ext _R ^1 (M, Y_{i-1}) =0
$$
for $i>0$, as desired.
\qed
\end{pf}

Now we are in the position to prove the main result of this section.

\begin{thm}\label{gcontfinc}
Let $R$ be a generically Gorenstein Cohen-Macaulay local ring with canonical module.
Suppose that $\gt\cap\C$ is covariantly finite in $\C$.
Then $\G$ is contravariantly finite in $\C$, hence in $\mod\,R$.
\end{thm}

\begin{pf}
It follows from Proposition \ref{gtthick} and \cite[Proposition (1.3)]{Yoshino2} that both $\gt$ and $\C$ are closed under extensions.
Hence $\gt\cap\C$ is also closed under extensions, and we see from Corollary \ref{tycont} that ${}^{\sf L} (\gt\cap\C)\cap\C$ is contravariantly finite in $\C$.
Using Lemma \ref{cateq}, Corollary \ref{gengorcm} and Proposition \ref{gcm}, we get
$$
{}^{\sf L}(\gt\cap\C)\cap\C = {}^{\bot}(\gt\cap\C)\cap\C = \tgt\cap\C = \G\cap\C = \G.
$$
Therefore $\G$ is contravariantly finite in $\C$.
Since $\C$ is contravariantly finite in $\mod\,R$ by \cite[Corollary (4.20)]{Yoshino2}, we see that $\G$ is contravariantly finite in $\mod\,R$.
\qed
\end{pf}

The author proved that if $R$ is a non-Gorenstein local ring with $\depth\,R\leq 2$ and $\G\neq\F$, then there exists a module which does not admit a right $\G$-approximation.
For the details, see \cite{Takahashi1}, \cite{Takahashi2} and \cite{Takahashi3}.

\begin{lem}\label{mymain}
Let $(R, \m, k)$ be a non-Gorenstein local ring with $\G\neq\F$.
\begin{enumerate}
\item[{\rm (1)}]
If $\depth\,R=0$, then $k$ does not have a right $\G$-approximation.
\item[{\rm (2)}]
If $\depth\,R=1$, then $\m$ does not have a right $\G$-approximation.
\item[{\rm (3)}]
If $\depth\,R=2$ and $0 \to R \to E \to \m \to 0$ is a nonsplit exact sequence, then $E$ does not have a right $\G$-approximation.
\end{enumerate}
\end{lem}

Using this lemma, as a corollary of the above theorem we get the following peculiar result.

\begin{cor}\label{infinity}
Let $R$ be a generically Gorenstein Cohen-Macaulay local ring with canonical module.
Suppose that $R$ is non-Gorenstein, $\dim R\leq 2$ and $\G\neq\F$.
Then there exist infinitely many nonisomorphic indecomposable G-projective modules $M$ and infinitely many nonisomorphic indecomposable maximal Cohen-Macaulay modules $N$ such that $\Ext _R ^i (M, N) = 0$ for all $i>0$.
\end{cor}

\begin{pf}
It is seen from Lemma \ref{mymain} that $\G$ is not contravariantly finite in $\mod\,R$.
Hence $\gt\cap\C$ is not covariantly finite in $\C$ by Theorem \ref{gcontfinc}.
Therefore both $\G$ and $\gt\cap\C$ contain infinitely many nonisomorphic indecomposable $R$-modules by \cite[Proposition 4.2]{AS}.
This proves the corollary.
\qed
\end{pf}

There actually exists a local ring $R$ satisfying the assumptions of the above corollary, as follows:

\begin{ex}
Let
$$
R=k[[X,Y,Z,W]]/(X^2,Y^2-YW,YZ-YW,Z^2-YW),
$$
where $k$ is a field.
Denote by $x,y,z,w$ the residue classes of $X,Y,Z,W$ in $R$ respectively.
Then $R$ is a one-dimensional complete Cohen-Macaulay non-Gorenstein local ring with parameter $w$, and the minimal primes of $R$ are $\p = (x,y,z)$, $\q =(x,y-w,z-w)$.
It is easy to observe that the local rings $R_{\p}$ and $R_{\q}$ are complete intersections, hence Gorenstein rings.
Therefore $R$ is generically Gorenstein.
Since one has $(0: x)=(x)$, the $R$-module $R/(x)$ has a complete resolution
$$
\cdots \overset{x}{\to} R \overset{x}{\to} R \overset{x}{\to} \cdots.
$$
Hence $R/(x)$ is a nonfree G-projective $R$-module by Proposition \ref{g}(1), therefore one has $\G\neq\F$.
Thus the local ring $R$ satisfies the assumptions of Corollary \ref{infinity}.
\end{ex}

\section{The structure of $\rap\,\G$}

In this section, we will mainly study the structure of the modules of which there exist right $\G$-approximations.
We shall analyze the subcategory of $\mod\,R$ consisting of all such modules.

\begin{defn}
We define $\rap\,\G$ (resp. $\lap\,\G$) as the subcategory of $\mod\,R$ consisting of all $R$-modules that have right (resp. left) $\G$-approximations.
\end{defn}

Note that $\rap\,\G =\mod\,R$ (resp. $\lap\,\G =\mod\,R$) if and only if $\G$ is contravariantly finite (resp. covariantly finite) in $\mod\,R$.

We begin with giving a common property of $\rap\,\G$ and $\lap\,\G$.

\begin{prop}\label{rapsum}
Both $\rap\,\G$ and $\lap\,\G$ are subcategories of $\mod\,R$ containing $\G$ which are closed under finite direct sums and direct summands.
\end{prop}

\begin{pf}
We show only the assertion concerning $\rap\,\G$.
(The assertion concerning $\lap\,\G$ can be shown similarly.)
For any object $X$ of $\G$, the identity map $X \to X$ is a right $\G$-approximation of $X$.
Hence $\G\subseteq\rap\,\G$.
Let $M_1, M_2\in\mod\,R$.
Suppose that $f_1:X_1\to M_1$ and $f_2:X_2\to M_2$ are right $\G$-approximations of $M_1$ and $M_2$, respectively.
Then we easily see that the homomorphism
$$
\left(
\begin{smallmatrix}
f_1 & 0 \\
0 & f_2
\end{smallmatrix}
\right)
: X_1\oplus X_2 \to M_1\oplus M_2
$$
be a right $\G$-approximation of $M_1\oplus M_2$.
Hence $\rap\,\G$ is closed under finite direct sums.
On the other hand, suppose that $f : X \to M_1\oplus M_2$ is a right $\G$-approximation of $M_1\oplus M_2$.
Write $f = \binom{f_1}{f_2}$ along the decomposition.
Then we easily see that $f_1:X\to M_1$ and $f_2:X\to M_2$ are right $\G$-approximations of $M_1$ and $M_2$, respectively.
Hence $\rap\,\G$ is closed under direct summands.
\qed
\end{pf}

From now on, we set our sight on $\rap\,\G$.
It possesses the following properties.

\begin{prop}\label{raplem}
\begin{enumerate}
\item[{\rm (1)}]
$\rap\,\G$ contains $\gt$.
\item[{\rm (2)}]
$\rap\,\G$ is a resolving subcategory of $\mod\,R$.
\item[{\rm (3)}]
An $R$-module $M$ belongs to $\rap\,\G$ if and only if so does $\Omega M$.
\end{enumerate}
\end{prop}

\begin{pf}
(1) Let $M$ be an $R$-module in $\gt$.
Then we have an exact sequence
$$
0 \to \Omega M \to F \overset{\varepsilon}\to M \to 0
$$
where $F$ is a free $R$-module, hence $F$ belongs to $\G$.
The module $\Omega M$ belongs to $\gt$ because $\gt$ is resolving by Proposition \ref{gtthick}.
Therefore we see from the above exact sequence that the homomorphism $\varepsilon$ is a right $\G$-approximation of $M$, and thus $M$ is in $\rap\,\G$.

(2) It follows from Proposition \ref{gra} and \cite[Proposition 3.7(a)]{AR} that $\rap\,\G$ is closed under extensions.
According to Proposition \ref{rapsum} and \cite[Lemma 3.2(2)]{Yoshino}, we have only to show that $\rap\,\G$ is closed under syzygies.
Fix $M\in\rap\,\G$.
We have an exact sequence $0 \to Y \to X \to M \to 0$ with $X\in\G$ and $Y\in\gt$ by Proposition \ref{resappr}(2).
Taking the syzygies, we get an exact sequence
$$
0 \to \Omega Y \to \Omega X\oplus F \overset{\phi}{\to} \Omega M \to 0,
$$
where $F$ is free.
Since both $\G$ and $\gt$ are resolving by Propositions \ref{gra} and \ref{gtthick}, it follows that $\Omega X\oplus F$ and $\Omega Y$ belong to $\G$ and $\gt$, respectively.
Hence it is seen from the above exact sequence that $\phi$ is a right $\G$-approximation, which implies that $\Omega M$ belongs to $\rap\,\G$.

(3) The ``only if'' part was proved in (2).
Let $M$ be an $R$-module such that $\Omega M\in\rap\,\G$.
We want to show that $M\in\rap\,\G$.
According to Proposition \ref{resappr}(2), there is an exact sequence
$$
0 \to Y \to X \overset{f}{\to} \Omega M \to 0,
$$
where $X\in\G$ and $Y\in\gt$.
Proposition \ref{envelope3} yields the following diagram with exact rows:
$$
\begin{CD}
0 @>>> X @>{\rho}>> R^m @>>> \Omega ^{-1}X @>>> 0 \\
@. @V{f}VV \\
0 @>>> \Omega M @>>> R^n @>>> M @>>> 0
\end{CD}
$$
Note that $\rho$ is a left $\F$-approximation.
Hence the homomorphism $f$ can be lifted as follows:
$$
\begin{CD}
0 @>>> X @>{\rho}>> R^m @>>> \Omega ^{-1}X @>>> 0 \\
@. @V{f}VV @VVV @VVV \\
0 @>>> \Omega M @>>> R^n @>>> M @>>> 0
\end{CD}
$$
Adding some copies of $R$ to the first row in this diagram, we obtain a commutative diagram
$$
\begin{CD}
0 @>>> X @>>> R^{m+l} @>>> \Omega ^{-1}X\oplus R^l @>>> 0 \\
@. @V{f}VV @VVV @VVV \\
0 @>>> \Omega M @>>> R^n @>>> M @>>> 0 \\
@. @VVV @VVV @VVV \\
@. 0 @. 0 @. 0
\end{CD}
$$
with exact rows and columns.
Taking the kernels of the vertical maps, we get the following commutative diagram with exact rows and columns:
$$
\begin{CD}
@. 0 @. 0 @. 0 \\
@. @VVV @VVV @VVV \\
0 @>>> Y @>>> R^{m+l-n} @>>> Y' @>>> 0 \\
@. @VVV @VVV @VVV \\
0 @>>> X @>>> R^{m+l} @>>> \Omega ^{-1}X\oplus R^l @>>> 0 \\
@. @V{f}VV @VVV @VVV \\
0 @>>> \Omega M @>>> R^n @>>> M @>>> 0 \\
@. @VVV @VVV @VVV \\
@. 0 @. 0 @. 0
\end{CD}
$$
Proposition \ref{g}(2) says that $\Omega ^{-1}X\oplus R^l$ is in $\G$.
On the other hand, since $Y$ belongs to $\gt$, we see by the definition of $\gt$ and the long exact sequence of Ext that $Y'$ also belongs to $\gt$.
Hence, it follows from the exact sequence $0 \to Y' \to \Omega ^{-1}X\oplus R^l \to M \to 0$ in the above diagram that $M$ is in $\rap\,\G$.
\qed
\end{pf}

\begin{rem}
As we observed in Corollary \ref{abg}(1), all the modules of finite G-dimension admit right $\G$-approximations.
At first sight, it seems that no module of infinite G-dimension admits a right $\G$-approximation.
However, it is not true.
In fact, let $R$ be a Cohen-Macaulay non-Gorenstein local ring with canonical module $\omega$.
Then $\omega$ has infinite G-dimension because any Cohen-Macaulay local ring whose canonical module has finite G-dimension is Gorenstein (cf. \cite[Corollary 5.7]{ATY}), but $\omega$ has a right $\G$-approximation because $\omega$ belongs to $\gt$ by Proposition \ref{gcm} and $\gt$ is contained in $\rap\,\G$ by Proposition \ref{raplem}(1).
\end{rem}

In relation to the above remark, the condition that no module of infinite G-dimension admits a right $\G$-approximation can be translated as follows.

\begin{prop}
The following are equivalent:
\begin{enumerate}
\item[{\rm (1)}]
$\rap\,\G = \widehat{\G}$;
\item[{\rm (2)}]
$\gt\subseteq\widehat{\G}$;
\item[{\rm (3)}]
$\gt =\widehat{\F}$.
\end{enumerate}
\end{prop}

\begin{pf}
(1) $\Rightarrow$ (2): By Proposition \ref{raplem}(1) we have $\gt\subseteq\rap\,\G = \widehat{\G}$.

(2) $\Rightarrow$ (3): This implication follows from Corollary \ref{abg}(2).

(3) $\Rightarrow$ (1): Corollary \ref{abg}(1) yields the inclusion relation $\widehat{\G}\subseteq\rap\,\G$.
Conversely, let $M\in\rap\,\G$.
Then we have an exact sequence $0 \to Y \to X \to M \to 0$ with $X\in\G$ and $Y\in\gt$ by Proposition \ref{resappr}(2).
Since $\gt = \widehat{\F}\subseteq\widehat{\G}$, both of the modules $X$ and $Y$ are of finite G-dimension.
Hence we see that $M$ is also of finite G-dimension.
\qed
\end{pf}

As an application of Proposition \ref{raplem}, one can make from Lemma \ref{mymain} a series of modules which do not have right $\G$-approximations.

\begin{cor}\label{depamtwo}
Let $R$ be a non-Gorenstein local ring with $\depth\,R\leq 2$ and $\G\neq\F$.
Then $\Omega ^i k$ does not have a right $\G$-approximation for every $i\geq 0$.
\end{cor}

\begin{pf}
Suppose that $\depth\,R=0$ (resp. $\depth\,R=1$).
Then we have $k\not\in\rap\,\G$ (resp. $\Omega k=\m\not\in\rap\,\G$) by Lemma \ref{mymain}, hence $\Omega ^i k\not\in\rap\,\G$ for any $i\geq 0$ by Proposition \ref{raplem}(3).
Suppose that $\depth\,R=2$.
Then since $\Ext _R ^1 (\Omega k , R)\cong\Ext _R ^2 (k, R)\neq 0$, there exists a nonsplit exact sequence
$$
0 \to R \to E \to \Omega  k \to 0
$$
of $R$-modules.
Lemma \ref{mymain} says that $E$ does not belong to $\rap\,\G$.
On the other hand, $R$ belongs to $\rap\,\G$ and $\rap\,\G$ is closed under extensions by Proposition \ref{raplem}(2).
It follows that $\Omega k$ does not belong to $\rap\,\G$, and neither does $\Omega ^i k$ for any $i\geq 0$.
\qed
\end{pf}

For each module having a right $\G$-approximation, one can make three exact sequences associated to the module.

\begin{lem}\label{rapggen}
Let $M$ be an $R$-module in $\rap\,\G$.
Then there exist three exact sequences
$$
\begin{cases}
0 \to Y \to X \to M \to 0,\\
0 \to M \to Y' \to X' \to 0,\\
0 \to X \to M\oplus F \to Y' \to 0
\end{cases}
$$
in $\mod\,R$, where $X, X'\in\G$, $Y, Y'\in\gt$ and $F$ is free.
\end{lem}

\begin{pf}
It follows from Proposition \ref{resappr}(2) that we have an exact sequence $0 \to Y \to X \to M \to 0$ with $X\in\G$ and $Y\in\gt$.
Noting that $\F$ is a cogenerator for $\G$ by Proposition \ref{gra} again, we get an exact sequence $0 \to X \to F \to X' \to 0$, where $F$ is a free $R$-module and $X'$ is in $\G$.
Thus we obtain the pushout diagram:
$$
\begin{CD}
@. @. 0 @. 0 \\
@. @. @VVV @VVV \\
0 @>>> Y @>>> X @>>> M @>>> 0 \\
@. @| @VVV @VVV \\
0 @>>> Y @>>> F @>>> Y' @>>> 0 \\
@. @. @VVV @VVV \\
@. @. X' @= X' \\
@. @. @VVV @VVV \\
@. @. 0 @. 0
\end{CD}
$$
Since $\gt$ is thick and contains $R$ by Proposition \ref{gtthick}, the module $Y'$ belongs to $\gt$.
We have the following pullback diagram:
$$
\begin{CD}
@. @. 0 @. 0 \\
@. @. @VVV @VVV \\
@. @. X @= X \\
@. @. @VVV @VVV \\
0 @>>> M @>>> P @>>> F @>>> 0 \\
@. @| @VVV @VVV \\
0 @>>> M @>>> Y' @>>> X' @>>> 0 \\
@. @. @VVV @VVV \\
@. @. 0 @. 0
\end{CD}
$$
The second row in the above diagram splits as $F$ is free, and thus we obtain an exact sequence $0 \to X \to M\oplus F \to Y' \to 0$.
\qed
\end{pf}

We have reached the stage to prove our main theorem in this section.
The structure of the subcategory $\rap\,\G$ is as follows.

\begin{thm}\label{rapgthick}
\begin{enumerate}
\item[{\rm (1)}]
$\rap\,\G$ is the smallest subcategory of $\mod\,R$ containing $\G$ and $\gt$ and closed under direct summands and extensions.
\item[{\rm (2)}]
$\rap\,\G$ is a thick subcategory of $\mod\,R$.
\end{enumerate}
\end{thm}

\begin{pf}
(1) Propositions \ref{rapsum} and \ref{raplem}(1) imply that both $\G$ and $\gt$ are contained in $\rap\,\G$.
Since $\rap\,\G$ is resolving by Proposition \ref{raplem}(2), $\rap\,\G$ is closed under direct summands and extensions.
On the other hand, letting $M\in\rap\,\G$, we have an exact sequence
$$
0 \to X \to M\oplus F \to Y \to 0
$$
in $\mod\,R$ with $X\in\G$ and $Y\in\gt$ by virtue of Lemma \ref{rapggen}.
Thus the assertion is proved.

(2) By Proposition \ref{raplem}(2), we have only to show that $\rap\,\G$ is closed under cokernels of monomorphisms.
Let $0 \to L \to M \to N \to 0$ be an exact sequence of $R$-modules with $L, M\in\rap\,\G$.
Taking the syzygy of $N$, one gets an exact sequence $0 \to \Omega N \to R^n \to N \to 0$.
From these exact sequences, one obtains the following pullback diagram:
$$
\begin{CD}
@. @. 0 @. 0 \\
@. @. @VVV @VVV \\
@. @. \Omega N @= \Omega N \\
@. @. @VVV @VVV \\
0 @>>> L @>>> P @>>> R^n @>>> 0 \\
@. @| @VVV @VVV \\
0 @>>> L @>>> M @>>> N @>>> 0 \\
@. @. @VVV @VVV \\
@. @. 0 @. 0
\end{CD}
$$
Since the middle row in the diagram splits, we get an exact sequence
$$
0 \to \Omega N \to L\oplus R^n \to M \to 0.
$$
Since $\rap\,\G$ is closed under finite direct sums and contains $R$ by Proposition \ref{raplem}(2), $L\oplus R^n$ belongs to $\rap\,\G$, and so does $\Omega N$ because $\rap\,\G$ is closed under kernels of epimorphisms by Proposition \ref{raplem}(2) again.
Finally, using Proposition \ref{raplem}(3), we conclude that $N$ belongs to $\rap\,\G$, as desired.
\qed
\end{pf}

\begin{rem}
Using the above theorem, one can give another proof of the first statement of Corollary \ref{abg}(1):

Let $M$ be an $R$-module of finite G-dimension.
Then, by definition, we have an exact sequence
$$
0 \to X_n \to X_{n-1} \to \cdots \to X_1 \to X_0 \to M \to 0
$$
of $R$-modules with $X_i\in\G$ for $0\leq i\leq n$.
Firstly, decompose this exact sequence into short exact sequences.
Secondly, note from Proposition \ref{rapsum} that $\rap\,\G$ contains $\G$, and from Theorem \ref{rapgthick}(2) that $\rap\,\G$ is closed under cokernels of monomorphisms.
Then one sees that $M$ belongs to $\rap\,\G$.
\end{rem}

Let $R$ be an Artinian ring.
A subcategory $\X$ of $\mod\,R$ is called {\it coresolving} if the following three conditions hold:
\begin{enumerate}
\item[{\rm (1)}]
$\X$ contains all injective $R$-modules,
\item[{\rm (2)}]
$\X$ is closed under extensions,
\item[{\rm (3)}]
$\X$ is closed under cokernels of monomorphisms.
\end{enumerate}
We end this section by remarking that the subcategory $\rap\,\G$ is not only resolving but also coresolving over an Artinian ring $R$.

\begin{cor}
Suppose that $R$ is Artinian.
Then $\rap\,\G$ is a coresolving subcategory of $\mod\,R$.
\end{cor}

\begin{pf}
According to Theorem \ref{rapgthick}(2), it is enough to prove that $\rap\,\G$ contains all injective $R$-modules.
However, it is obvious because any injective $R$-module belongs to $\gt$ and $\gt$ is contained in $\rap\,\G$ by Proposition \ref{raplem}(1).
\qed
\end{pf}

\section{A characterization in terms of Tate cohomologies}

In this section, we will make a characterization of the subcategory $\rap\,\G$ in terms of Tate cohomologies.
To be more concrete, we shall give a criterion for an $R$-module to admit a right $\G$-approximation by the vanishing of certain Tate cohomology modules.
Before stating the definition of a Tate cohomology module, we introduce the notion of the stable category of a given additive category, and give several related results.

For a subcategory $\X$ of $\mod\,R$, we denote by $\underline{\X}$ the {\it stable category} of $\X$, namely, the objects of $\underline{\X}$ are the same as those of $\X$, and for objects $M, N$ of $\X$, the set of morphisms from $M$ to $N$ is defined by
$$
\lhom _R (M, N):=\Hom _R (M, N)/\P _R(M, N),
$$
where $\P _R(M, N)$ is the $R$-submodule of $\Hom _R (M, N)$ consisting of all homomorphisms from $M$ to $N$ factoring through some free $R$-module.
We denote by $\underline{f}$ the residue class of $f\in\Hom _R (M, N)$ in $\lhom _R (M, N)$.

Let $f: M \to N$ be a homomorphism of $R$-modules.
Then we see that there are commutative diagrams
$$
\begin{CD}
0 @>>> \Omega M @>>> P_0 @>>> M @>>> 0 \\
@. @V{g}VV @VVV @V{f}VV \\
0 @>>> \Omega N @>>> Q_0 @>>> N @>>> 0
\end{CD}
$$
$$
\begin{CD}
M @>>> P^0 @>>> \Omega ^{-1}M @>>> 0 \\
@V{f}VV @VVV @V{h}VV \\
N @>>> Q^0 @>>> \Omega ^{-1}N @>>> 0
\end{CD}
$$
with exact rows, where $P_0, Q_0, P^0, Q^0$ are free $R$-modules.
On the other hand, we have a commutative diagram
$$
\begin{CD}
P_1 @>>> P_0 @>>> M @>>> 0 \\
@VVV @VVV @V{f}VV \\
Q_1 @>>> Q_0 @>>> N @>>> 0
\end{CD}
$$
where the rows are minimal free presentations of $M$ and $N$, respectively.
Dualizing this diagram by $R$ gives the following commutative diagram:
$$
\begin{CD}
0 @>>> N^{\ast} @>>> Q_0^{\ast} @>>> Q_1^{\ast} @>>> \tr N @>>> 0 \\
@. @V{f^{\ast}}VV @VVV @VVV @V{l}VV \\
0 @>>> M^{\ast} @>>> P_0^{\ast} @>>> P_1^{\ast} @>>> \tr M @>>> 0
\end{CD}
$$
It is easy to check that the homomorphisms $g, h, l$ are uniquely determined up to homomorphism factoring through some free $R$-module, and that if $f$ factors through a free $R$-module, then one can choose the zero maps as $g, h, l$.
Thus the homomorphisms
$$
\begin{cases}
\lhom _R (M, N) \longrightarrow \lhom _R (\Omega M, \Omega N),\\
\lhom _R (M, N) \longrightarrow \lhom _R (\Omega ^{-1}M, \Omega ^{-1}N),\\
\lhom _R (M, N) \longrightarrow \lhom _R (\tr N, \tr M)
\end{cases}
$$
given by $\underline{f}\mapsto\underline{g}$, $\underline{f}\mapsto\underline{h}$, $\underline{f}\mapsto\underline{l}$ respectively, are well-defined.
We should note that the third homomorphism is an isomorphism since any $R$-module $L$ is isomorphic to $\tr (\tr L)$ up to free summand.

The above observation means that $\Omega, \Omega ^{-1}$ define functors from $\underline{\mod\,R}$ to itself, and $\tr$ defines a functor from $(\underline{\mod\,R})^{\rm op}$ to $\underline{\mod\,R}$ giving an equivalence of categories.

The functors $\Omega, \Omega ^{-1}$ behave well on the stable category of $\G$, as follows.
One can easily prove this proposition by using \cite[Proposition (2.46)]{ABr}.

\begin{prop}\label{omegainverse}
For G-projective $R$-modules $M, N$, the homomorphisms
$$
\begin{cases}
\lhom _R (M, N) \longrightarrow \lhom _R (\Omega M, \Omega N) \\
\lhom _R (M, N) \longrightarrow \lhom _R (\Omega ^{-1}M, \Omega ^{-1}N)
\end{cases}
$$
defined by $\Omega$ and $\Omega ^{-1}$ respectively, are isomorphisms.
One of the homomorphisms is the inverse map of the other.

In other words, $\Omega$ defines an isomorphic functor from $\lg$ to itself with the inverse functor $\Omega ^{-1}$.
\end{prop}

The $R$-module $\lhom _R (X, M)$ can be represented by Ext modules if $X$ belongs to $\G$.

\begin{lem}\label{homext}
Let $M$ be an $R$-module, and $X$ a G-projective $R$-module.
Then
$$
\lhom _R (X, M) \cong \Ext _R ^1 (X, \Omega M) \cong \Ext _R ^1 (\Omega ^{-1}X, M).
$$
\end{lem}

\begin{pf}
Take the first syzygy of $M$; one has an exact sequence
$$
0 \to \Omega M \to F \overset{\pi}{\to} M \to 0
$$
where $F$ is a free $R$-module.
Applying the functor $\Hom _R (X, -)$ to this exact sequence, one gets an exact sequence
$$
\Hom _R (X, F) \overset{\rho}{\to} \Hom _R (X, M) \to \Ext _R ^1 (X, \Omega M) \to 0.
$$
Note from Proposition \ref{freecover}(1) that $\pi$ is a right $\F$-approximation of $M$.
It is easily seen that the image of the map $\rho =\Hom _R (X, \pi )$ coincides with $\P _R (X, M)$.
Thus an isomorphism
$$
\lhom _R (X, M) \cong \Ext _R ^1 (X, \Omega M).
$$
is obtained.

As for the other isomorphism, by Proposition \ref{envelope3}, there is an exact sequence
$$
0 \to X \overset{\theta}{\to} F' \to \Omega ^{-1}X \to 0,
$$
where $\theta$ is a left $\F$-approximation.
Dualizing this sequence by $M$ gives an exact sequence
$$
\Hom _R (F', M) \overset{\kappa}{\to} \Hom _R (X, M) \to \Ext _R ^1 (\Omega ^{-1}X, M) \to 0,
$$
and the image of $\kappa =\Hom _R (\theta, M)$ coincides with $\P _R (X, M)$.
Thus one gets an isomorphism
$$
\lhom _R (X, M) \cong \Ext _R ^1 (\Omega ^{-1}X, M).
$$
\qed
\end{pf}

Let $X\in\G$ and $M\in\mod\,R$.
For each $i\in\Z$, we define the {\it $i$th Tate cohomology module} by
$$
\cext _R ^i (X, M) = \lhom _R (\Omega ^i X, M).
$$
Note that we have $\cext _R ^0 (X, M) = \lhom _R (X, M)$.

Let us study several basic properties of Tate cohomology modules.

\begin{prop}\label{tate}
\begin{enumerate}
\item[{\rm (1)}]
Let $M$ be an $R$-module and $X$ a G-projective $R$-module. 
\begin{enumerate}
\item[{\rm (i)}]
For $i, n\in\Z$ with $n>0$, one has
$$
\cext _R ^i (X, M) \cong \Ext _R ^n (\Omega ^{i-n} X, M).
$$
In particular, $\cext _R ^n (X, M) \cong \Ext _R ^n (X, M)$.
\item[{\rm (ii)}]
Let $F_{\bullet}$ be a complete resolution of $X$.
Then
$$
\cext _R ^i (X, M) \cong \H ^i (\Hom _R (F_{\bullet}, M))
$$
for $i\in\Z$.
\end{enumerate}
\item[{\rm (2)}]
\begin{enumerate}
\item[{\rm (i)}]
Let $X$ be a G-projective $R$-module, and let $0 \to M' \to M \to M'' \to 0$ be an exact sequence of $R$-modules.
Then there is a long exact sequence
$$
\begin{CD}
\cdots @>>> \cext _R ^i (X, M') @>>> \cext _R ^i (X, M) @>>> \cext _R ^i (X, M'') \\
@>>> \cext _R ^{i+1} (X, M') @>>> \cdots \quad (i\in\Z).
\end{CD}
$$
\item[{\rm (ii)}]
Let $M$ be an $R$-module, and let $0 \to X' \to X \to X'' \to 0$ be an exact sequence of G-projective $R$-modules.
Then there is a long exact sequence
$$
\begin{CD}
\cdots @>>> \cext _R ^i (X'', M) @>>> \cext _R ^i (X, M) @>>> \cext _R ^i (X', M) \\
@>>> \cext _R ^{i+1} (X'', M) @>>> \cdots \quad (i\in\Z).
\end{CD}
$$
\end{enumerate}
\end{enumerate}
\end{prop}

\begin{pf}
(1)(i) We have $\cext ^i (X, M) = \lhom (\Omega ^i X, M)$ and $\Ext ^n (\Omega ^{i-n}X, M)\cong\Ext ^1 (\Omega ^{i-1}X, M)$.
It is seen from Proposition \ref{omegainverse} that $\Omega ^{i-1} X$ is isomorphic to $\Omega ^{-1}\Omega ^i X$ up to free summand.
Applying Lemma \ref{homext} to the G-projective module $\Omega ^i X$ (cf. Proposition \ref{g}(2)), we get $\lhom (\Omega ^i X, M)\cong\Ext ^1 (\Omega ^{-1}\Omega ^iX, M)\cong\Ext ^1 (\Omega ^{i-1}X, M)$.
Thus we obtain an isomorphism $\cext ^i (X, M)\cong\Ext ^n (\Omega ^{i-n} X, M)$.

(ii) Proposition \ref{g}(1) guarantees that $X$ has a complete resolution.
By the assertion (i), we have $\cext ^i (X, M)\cong\Ext ^1 (\Omega ^{i-1} X, M)$.
We see from Proposition \ref{compresol} that the image of the $(i-1)$th differential map of $F_{\bullet}$ is isomorphic to $\Omega ^{i-1}X$ up to free summand.
Noting this, we obtain an isomorphism $\H ^i (\Hom (F_{\bullet}, M)) \cong \Ext ^1 (\Omega ^{i-1}X, M)$.

(2)(i) Applying the functor $\Hom (\Omega ^{i-1}X, -)$ to the given short exact sequence, we get a long exact sequence:
\begin{equation*}
\begin{split}
\cdots & \to \Ext ^1 (\Omega ^{i-1}X, M') \to \Ext ^1 (\Omega ^{i-1}X, M) \to \Ext ^1 (\Omega ^{i-1}X, M'') \\
& \to \Ext ^2 (\Omega ^{i-1}X, M') \to \cdots.
\end{split}
\end{equation*}
Using the statement (1), we see that this gives the long exact sequence which we want.

(ii) For each integer $i$, there is an exact sequence of this form:
$$
0 \to \Omega ^{i-1}X' \to \Omega ^{i-1}X\oplus R^m \to \Omega ^{i-1}X'' \to 0.
$$
Dualizing this sequence by $M$, one gets a long exact sequence:
\begin{equation*}
\begin{split}
\cdots & \to \Ext ^1 (\Omega ^{i-1}X'', M) \to \Ext ^1 (\Omega ^{i-1}X, M) \to \Ext ^1 (\Omega ^{i-1}X', M) \\
& \to \Ext ^2 (\Omega ^{i-1}X'', M) \to \cdots.
\end{split}
\end{equation*}
It follows from (1) this can be identified with the exact sequence in the assertion.
\qed
\end{pf}

\begin{rem}
Let $M$ be an $R$-module and $X$ a G-projective $R$-module.
Avramov and Martsinkovsky \cite{AM} defines the $i$th Tate cohomology module by $\cext _R ^i (X, M) = \H ^i (\Hom _R (F_{\bullet}, M))$, where $F_{\bullet}$ is a complete resolution of $X$.
Proposition \ref{tate}(1)(ii) says that their definition is the same as ours.
\end{rem}

The main result of this section is the following theorem.

\begin{thm}\label{sensegeneral}
The following are equivalent for an $R$-module $M$:
\begin{enumerate}
\item[{\rm (1)}]
$M\in\rap\,\G$;
\item[{\rm (2)}]
$\cext _R ^i (-, M)|_{\lg}$ is a finitely generated right $\lg$-module for every $i\in\Z$;
\item[{\rm (2')}]
$\cext _R ^i (-, M)|_{\lg}$ is a finitely generated right $\lg$-module for some $i\in\Z$;
\item[{\rm (3)}]
$\cext _R ^i (-, M)|_{\lg}$ is a finitely presented right $\lg$-module for every $i\in\Z$;
\item[{\rm (3')}]
$\cext _R ^i (-, M)|_{\lg}$ is a finitely presented right $\lg$-module for some $i\in\Z$;
\item[{\rm (4)}]
$\cext _R ^i (-, M)|_{\lg}$ is a projective object of $\mod\,\lg$ for every $i\in\Z$;
\item[{\rm (4')}]
$\cext _R ^i (-, M)|_{\lg}$ is a projective object of $\mod\,\lg$ for some $i\in\Z$.
\end{enumerate}
\end{thm}

\begin{pf}
The implications (4) $\Rightarrow$ (3) $\Rightarrow$ (2), (4') $\Rightarrow$ (3') $\Rightarrow$ (2'), (4) $\Rightarrow$ (4'), (3) $\Rightarrow$ (3') and (2) $\Rightarrow$ (2') are obvious.
It is enough to show the implications (2') $\Rightarrow$ (1) $\Rightarrow$ (4).

(2') $\Rightarrow$ (1): For some integer $i$, there is an epimorphism
$$
\phi : \lhom (-, X)|_{\lg} \to \cext ^i (-, M)|_{\lg},
$$
where $X$ is a G-projective $R$-module.
We have a surjective homomorphism
$$
\phi (X): \lhom (X, X) \to \cext ^i(X, M)=\lhom (\Omega ^i X, M),
$$
and $\phi (X)(\underline{\id _X})=\underline{f_0}$ for some $f_0\in\Hom (\Omega ^i X, M)$.
Let $f_1: F \to M$ be a surjective homomorphism from a free $R$-module $F$, and set $f=(f_0, f_1): \Omega ^iX\oplus F \to M$.
Note then that $f$ is a surjective homomorphism satisfying $\underline{f}=\underline{f_0}$.

Let us show that $f$ is a right $\G$-approximation of $M$.
Take a homomorphism $f': X' \to M$ such that $X'$ is a G-projective $R$-module.
Note that $\cext ^i (\Omega ^{-i}X', M)$ can be identified with $\lhom (X', M)$ (cf. Proposition \ref{omegainverse}).
The surjectivity of $\phi (\Omega ^{-i}X')$ implies that there exists $g_0\in\Hom (\Omega ^{-i}X', X)$ such that $\underline{f'}=\phi (\Omega ^{-i}X')(\underline{g_0})$.
On the other hand, since $\phi$ is a natural transformation, we have the following commutative diagram:
$$
\begin{CD}
\lhom (X,X) @>{\phi (X)}>> \cext ^i(X,M) @= \lhom (\Omega ^i X, M) \\
@V{\lhom (\underline{g_0}, X)}VV @V{\cext ^i (\underline{g_0}, X)}VV @V{\lhom (\Omega ^i \underline{g_0}, M)}VV \\
\lhom (\Omega ^{-i}X', X) @>{\phi (\Omega ^{-i}X')}>> \cext ^i(\Omega ^{-i}X', M) @= \lhom (X', M)
\end{CD}
$$
The commutativity of this diagram yields $\underline{f'}=\underline{f}\cdot\Omega ^i\underline{g_0}$.
We can write $\Omega ^i \underline{g_0} = \underline{g}$ for some $g\in\Hom (X', \Omega ^i X\oplus F)$, and get $\underline{f'}=\underline{fg}$.
This means that the homomorphism $f'-fg$ factors through some free $R$-module $F'$; there exist $\alpha\in\Hom (X',F')$ and $\beta\in\Hom (F',M)$ such that $f'-fg=\beta\alpha$.
Noting that $f$ is a surjective homomorphism, we see that there exists $\gamma\in\Hom (F',\Omega ^iX\oplus F)$ satisfying $\beta =f\gamma$, and hence we have $f'=f(g+\gamma\alpha )$.
Thus the homomorphism $g$ factors through $f$, and we conclude that $f$ is a right $\G$-approximation of $M$.

(1) $\Rightarrow$ (4): By Proposition \ref{resappr}(2), there is an exact sequence $0 \to Y \to X \to M \to 0$ with $X\in\G$ and $Y\in\gt$.
Fix $X'\in\G$.
One has $\cext ^i (X', Y) = \lhom (\Omega ^i X', Y) \cong \Ext ^1 (\Omega ^i X', \Omega Y)$ by Lemma \ref{homext}, and $\Ext ^1 (\Omega ^i X', \Omega Y)=0$ because $\Omega ^iX'\in\G$ by Proposition \ref{g}(2) and $\Omega Y\in\gt$ by Proposition \ref{gtthick}.
Thus $\cext ^i(-, Y)|_{\lg}=0$, hence $\cext ^i (-, X)|_{\lg}\cong\cext ^i (-, M)|_{\lg}$ for any $i\in\Z$ by Proposition \ref{tate}(2)(i).
Since $\cext ^i (-, X)|_{\lg}=\lhom (\Omega ^i(-), X)|_{\lg} \cong \lhom (-, \Omega ^{-i}X)|_{\lg}$ by Proposition \ref{omegainverse}, the functor $\cext ^i (-, M)|_{\lg}$ is a projective object of $\mod\,\lg$.
\qed
\end{pf}

\begin{rem}\cite[Remark 2.6]{Yoshino}
Let $\X$ be a subcategory of $\mod\,R$, and let $\iota : \mod\,\underline{\X} \to \mod\,\X$ be the functor induced by the natural functor $\X\to\underline{\X}$.
Then $\iota$ gives an equivalence of categories between $\mod\,\underline{\X}$ and the full subcategory of $\mod\,\X$ consisting of all objects $F$ satisfying $F (R) =0$.
Thus, for example, one can identify the right $\lg$-module $\cext _R ^i (-, M)|_{\lg}$ with the right $\G$-module $\cext _R ^i (-, M)|_{\G}$.
\end{rem}

As an immediate corollary of the above theorem, we obtain a criterion for an $R$-module to have right $\G$-approximation in terms of $\lhom$.

\begin{cor}\label{sense}
The following are equivalent for an $R$-module $M$:
\begin{enumerate}
\item[{\rm (1)}]
$M\in\rap\,\G$;
\item[{\rm (2)}]
$\lhom _R (-, M)|_{\lg}$ is a finitely generated right $\lg$-module;
\item[{\rm (3)}]
$\lhom _R (-, M)|_{\lg}$ is a finitely presented right $\lg$-module;
\item[{\rm (4)}]
$\lhom _R (-, M)|_{\lg}$ is a projective object of $\mod\,\lg$.
\end{enumerate}
\end{cor}

\section{$\G$-approximations over reduced rings}

In this section, we will observe $\G$-approximations mainly over reduced rings.
Considering the relationships between $\rap\,\G$ and $\lap\,\G$, we shall give sufficient conditions for the covariant finiteness and contravariant finiteness of $\G$ in $\mod\,R$.

Let us start by showing the following easy lemma.

\begin{lem}\label{filt}
Let $\X$ be a subcategory of $\mod\,R$ which is closed under extensions, and let $M$ be an $R$-module.
Suppose that $R/\p$ belongs to $\X$ for any $\p\in\Supp _R\,M$.
Then $M$ belongs to $\X$.
\end{lem}

\begin{pf}
There is a filtration
$$
0 =M_0 \subseteq M_1 \subseteq \cdots \subseteq M_{n-1} \subseteq M_n =M
$$
of $R$-submodules of $M$ such that $M_i/M_{i-1}\cong R/\p _i$ for some $\p _i\in\Supp _R\, M$.
Decompose this filtration into short exact sequences.
Noting that $\X$ is closed under extensions and each $R/\p _i$ belongs to $\X$, we easily observe that $M$ belongs to $\X$.
\qed
\end{pf}

The following proposition will play a key role throughout this section.

\begin{prop}\label{positivegrade}
Let $R$ be a reduced ring, and let $\X$ be a subcategory of $\mod\,R$ containing $R$ which is closed under direct summands and extensions.
Suppose that any module $M$ with $M^{\ast}=0$ belongs to $\X$.
Then $\X = \mod\,R$.
\end{prop}

\begin{pf}
Fix $M\in\mod\,R$.
We want to show that $M$ belongs to $\X$.
By virtue of Lemma \ref{filt}, without loss of generality, we can assume $M=R/\p$ where $\p$ is a prime ideal of $R$.

For an ideal $I$ of $R$, we denote by $\lambda I$ the ideal $I+(0:I)$ of $R$.
Noting that $R$ has no nonzero nilpotents, one easily observes that $I\cap (0:I)=0$ for any ideal $I$.
Setting $J=(0:I)$, one has $(R/\lambda I)^{\ast}\cong (0:\lambda I)=(0:I+(0:I))=(0:I)\cap (0:(0:I))=J\cap (0:J)=0$.
The assumption of the proposition says that $R/\lambda I$ belongs to $\X$ for any ideal $I$ of $R$.

Since $\p\cap (0:\p )=0$, we have an exact sequence
$$
0 \to R \overset{f}{\to} R/\p \oplus R/(0:\p ) \overset{g}{\to} R/\lambda\p \to 0,
$$
where $f(a)=\binom{\overline{a}}{\overline{a}}$ and $g(\binom{\overline{x}}{\overline{y}})=\overline{x-y}$, and $R/\lambda\p$ belongs to $\X$.
Since $\X$ contains $R$ and is closed under extensions, the middle module $R/\p \oplus R/(0:\p )$ in the exact sequence also belongs to $\X$.
Since $\X$ is closed under direct summands, the $R$-module $R/\p$ also belongs to $\X$, as desired.
\qed
\end{pf}

\begin{rem}
Proposition \ref{positivegrade} does not necessarily hold unless the assumption that the base ring $R$ is reduced.
In fact, let $R=k[[x, y]]/(x^2)$ where $k$ is a field, and let $\X$ be the subcategory of $\mod\,R$ generated by $R$ and $k$ as a subcategory closed under summands and extensions.
Let $M$ be an $R$-module with $M^{\ast}=0$, equivalently, $\grade\,M>0$.
Then, noting that $R$ is a Cohen-Macaulay local ring of dimension one, we have $\grade\,M=\codim\,M=1-\Kdim\,M$.
Hence $\Kdim\,M=0$, in other words, $M$ has finite length.
Since $\X$ contains the $R$-module $k$, we see that $M$ belongs to $\X$ by Lemma \ref{filt}.
Thus $\X$ satisfies the assumptions of Proposition \ref{positivegrade}.

On the other hand, set $\p = xR$.
Note that $\p$ is a prime ideal of $R$.
Let us consider the subcategory
$$
\M := \{\, M\in\mod\, R\, |\, M_{\p}\text{ is }R_{\p}\text{-free}\,\}
$$
of $\mod\,R$.
It is obviously seen that $\M$ is closed under direct summands and extensions and contains both $R$ and $k$.
This means that $\X$ is contained $\M$.
But since $(R/\p)_{\p}=\kappa (\p)$ is not $R_{\p}$-free, the $R$-module $R/\p$ does not belong to $\M$.
This especially says that $\X$ does not coincides with $\mod\,R$.
\end{rem}

Proposition \ref{positivegrade} yields the following corollary.

\begin{cor}\label{last}
Let $R$ be a reduced ring and $\X$ a subcategory of $\mod\,R$ which is closed under summands and extensions.
Suppose that either of the following holds:
\begin{enumerate}
\item[{\rm (1)}]
$\X$ contains $\widehat{\F}$ and is closed under transpose;
\item[{\rm (2)}]
$\X$ contains $R$ and $R/\p\in\X$ for any $\p\in\Spec\,R-\Ass\,R$.
\end{enumerate}
Then $\X = \mod\,R$.
\end{cor}

\begin{pf}
(1) Let $M$ be an $R$-module with $M^{\ast}=0$.
Take a minimal free presentation $R^n \to R^m \to M \to 0$ of $M$, and dualizing this by $R$, we obtain an exact sequence
$$
0=M^{\ast} \to R^m \to R^n \to \tr M \to 0,
$$
which says that $\tr M$ has projective dimension at most one.
Hence $\tr M\in\widehat{\F}\subseteq\X$.
Also we have $F\in\widehat{\F}\subseteq\X$ for any free $R$-module $F$.
The $R$-module $M$ is isomorphic to $\tr (\tr M)$ up to free summand, and $\tr (\tr M)\in\X$ as $\X$ is closed under transpose.
Since $\X$ is closed under finite sums and summands, we see that $M$ belongs to $\X$.
Thus it follows from Proposition \ref{positivegrade} that $\X$ coincides with $\mod\,R$.

(2) Let $M$ be an $R$-module satisfying $M^{\ast}=0$, i.e., $\grade\, M>0$.
Then, since $\grade\,M=\inf\{\,\depth\,R_{\p}\,|\,\p\in\Supp\,M\,\}$ by \cite[Proposition 1.2.10(a)]{BH}, we have $\depth\,R_{\p}>0$, equivalently $\p\not\in\Ass\,R$, for every $\p\in\Supp\,M$.
The assumption of the corollary says that $R/\p$ is in $\X$ for every $\p\in\Supp\,M$.
Lemma \ref{filt} implies that $M$ is in $\X$.
Finally, Proposition \ref{positivegrade} shows that $\X$ coincides with $\mod\,R$.
\qed
\end{pf}

We have already observed in Corollaries \ref{depamtwo} and \ref{abg}(3) that in the case where $R$ has depth at most two, $\G$ is contravariantly finite in $\mod\,R$ if $k$ has a right $\G$-approximation.
As follows, when $R$ is one-dimensional and reduced, this fact can be shown more easily.

\begin{cor}\label{onedimred}
Let $(R, \m , k)$ be a one-dimensional reduced local ring.
\begin{enumerate}
\item[{\rm (1)}]
Let $\X$ be a subcategory of $\mod\,R$ containing $R$ and $k$ which is closed under summands and extensions.
Then $\X =\mod\,R$.
\item[{\rm (2)}]
If $k$ has a right $\G$-approximation, then $\G$ is contravariantly finite in $\mod\,R$.
\end{enumerate}
\end{cor}

\begin{pf}
(1) Since $R$ is reduced, $R$ satisfies Serre's condition $(S_1)$.
Hence $\Ass\,R=\Min\,R$.
As $R$ has dimension one, we have $\Spec\,R-\Ass\,R=\{\m\}$.
Thus the assertion follows from Corollary \ref{last}(2).

(2) According to Proposition \ref{raplem}(2), the subcategory $\rap\,\G$ of $\mod\,R$ contains $R$ and is closed under summands and extensions.
Hence the assertion follows from (1).
\qed
\end{pf}

Next, we shall investigate the relationship between $\rap\,\G$ and $\lap\,\G$; the transpose $\tr$ corresponds a module in one of them to a module in the other.

\begin{prop}\label{laprap}
\begin{enumerate}
\item[{\rm (1)}]
An $R$-module $M$ belongs to $\rap\,\G$ (resp. $\lap\,\G$) if and only if $\tr M$ belongs to $\lap\,\G$ (resp. $\rap\,\G$).
\item[{\rm (2)}]
The category $\rap\,\G$ is closed under transpose if and only if $\lap\,\G$ is contained in $\rap\,\G$.
\end{enumerate}
\end{prop}

To prove this proposition, we need a lemma:

\begin{lem}\label{leftleft}
The following are equivalent for an $R$-module $M$:
\begin{enumerate}
\item[{\rm (1)}]
$M$ has a left $\G$-approximation;
\item[{\rm (2)}]
$\Hom (M, -)|_{\G}$ is a finitely generated left $\G$-module;
\item[{\rm (3)}]
$\lhom (M, -)|_{\lg}$ is a finitely generated left $\lg$-module.
\end{enumerate}
\end{lem}

\begin{pf}
(1) $\Rightarrow$ (2): Let $\phi : M\to X$ be a left $\G$-approximation.
Then it is easily seen from definition that $\Hom (\phi , -) : \Hom (X,-)\to\Hom (M,-)$ is a surjective morphism.

(2) $\Rightarrow$ (3): There is an epimorphism $\Phi : \Hom (X, -)|_{\G}\to\Hom (M, -)|_{\G}$, and using the Yoneda Lemma (cf. \cite[III-2]{maclane}), one sees that this epimorphism is induced by some homomorphism $\phi : M \to X$.
Hence one gets an epimorphism $\lhom (\underline{\phi}, -)|_{\lg} : \lhom (X, -)|_{\lg}\to\lhom (M, -)|_{\lg}$.

(3) $\Rightarrow$ (1): There is an epimorphism $\Psi : \lhom (X, -)|_{\lg}\to\lhom (M, -)|_{\lg}$, and by the Yoneda lemma there is a homomorphism $\phi : M \to X$ such that $\Psi = \lhom (\underline{\phi}, -)|_{\lg}$.
Note from definition that $X$ is reflexive.
Hence, adding some free $R$-module to $X$, one may assume that the homomorphism $\phi ^{\ast} : X^{\ast} \to M^{\ast}$ is surjective.

Let $\phi ': M \to X'$ be a homomorphism of $R$-modules such that $X'$ is a G-projective $R$-module.
Since $\Psi (X')$ is surjective, there exists $f\in\Hom (X, X')$ such that $\underline{\phi '}=\underline{f\phi}$.
Hence the homomorphism $\phi '-f\phi$ factors through some free $R$-module; there exist a free module $F$ and homomorphisms $\alpha\in\Hom (M,F)$, $\beta\in\Hom (F,X')$ such that $\phi '-f\phi =\beta\alpha$.
Noting that $\phi ^{\ast}$ is surjective and that $\alpha ^{\ast}$ is a map from a free module, one has $\alpha ^{\ast}=\phi ^{\ast}g$ for some $g\in\Hom (F^{\ast}, X^{\ast})$, hence $\alpha ^{\ast\ast}=g^{\ast}\phi ^{\ast\ast}$.
Denote by $\lambda _L$ the natural homomorphism from $L$ to $L^{\ast\ast}$ for an $R$-module $L$.
Setting $h = (\lambda _F)^{-1}\cdot g^{\ast}\cdot\lambda _X$, one has $h\phi = (\lambda _F)^{-1}\cdot g^{\ast}\cdot\lambda _X\cdot\phi = (\lambda _F)^{-1}\cdot g^{\ast}\cdot\phi ^{\ast\ast}\cdot\lambda _M = (\lambda _F)^{-1}\cdot \alpha ^{\ast}\cdot\lambda _M = \alpha$, hence $\phi ' = f\phi +\beta\alpha = (f+\beta h)\phi$.
Thus $\phi '$ factors through $\phi$, and one concludes that $\phi$ is a left $\G$-approximation of $M$.
\qed
\end{pf}

Now we can prove Proposition \ref{laprap}.

\begin{lrpf}
First of all, note that an $R$-module $M$ is isomorphic to $\tr (\tr M)$ up to free summand, and that both $\rap\,\G$ and $\lap\,\G$ are subcategories of $\mod\,R$ containing $R$ closed under finite sums and summands by Proposition \ref{rapsum}.
Hence $M$ belongs to $\rap\,\G$ (resp. $\lap\,\G$) if and only if $\tr (\tr M)$ belongs to $\rap\,\G$ (resp. $\lap\,\G$).

(1) It is enough to show that an $R$-module $M$ belongs to $\rap\,\G$ if and only if $\tr M$ belongs to $\lap\,\G$.
The condition that $M$ belongs to $\rap\,\G$ is equivalent to the condition that $\lhom (-, M)|_{\lg}$ is a finitely generated right $\lg$-module by Corollary \ref{sense}.
Note that $\lhom (X, M)$ is isomorphic to $\lhom (\tr M, \tr X)$ for each G-projective $R$-module $X$, and that an $R$-module belongs to $\G$ if and only if so does its transpose by Proposition \ref{g}(2).
Hence the condition that $\lhom (-, M)|_{\lg}$ is a finitely generated right $\lg$-module is equivalent to the condition that $\lhom (\tr M, -)|_{\lg}$ is a finitely generated left $\lg$-module.
Thus the assertion follows from Lemma \ref{leftleft}.

(2) Assume that $\rap\,\G$ is closed under transpose.
Let $M\in\lap\,\G$.
Then $\tr M\in\rap\,\G$ by (1).
Hence $\tr (\tr M)\in\rap\,\G$ by the assumption, and therefore $M\in\rap\,\G$.
Thus $\lap\,\G$ is contained in $\rap\,\G$.
Conversely, if this is the case, then for $M\in\rap\,\G$ one has $\tr M\in\lap\,\G\subseteq\rap\,\G$ by (1) and the assumption.
Therefore $\rap\,\G$ is closed under transpose.
\qed
\end{lrpf}

Proposition \ref{laprap} together with Corollary \ref{abg}(3) yield the following:

\begin{cor}
If $R$ is Gorenstein, then $\G$ is functorially finite in $\mod\,R$.
\end{cor}

Now, let us achieve the main aim of this section; the following is the main result of this section.

\begin{thm}
Let $R$ be a reduced ring.
If $\lap\,\G\subseteq\rap\,\G$, then $\G$ is contravariantly finite in $\mod\,R$.
\end{thm}

\begin{pf}
Propositions \ref{rapsum} and \ref{laprap}(2) and the assumption imply that $\rap\,\G$ is closed under summands, extensions and transpose.
Corollary \ref{abg}(1) especially says that $\rap\,\G$ contains all $R$-modules of finite projective dimension.
Hence it follows from Corollary \ref{last}(1) that $\rap\,\G$ coincides with $\mod\,R$, which means that $\G$ is contravariantly finite in $\mod\,R$.
\qed
\end{pf}

We end the present paper with a sufficient condition for the covariant finiteness of $\G$ in $\mod\,R$ in the case where $R$ is a domain.

\begin{prop}
Let $R$ be an integral domain.
If $\lap\,\G$ is closed under extensions, then $\G$ is covariantly finite in $\mod\,R$.
\end{prop}

\begin{pf}
We want to prove that $\lap\,\G$ coincides with $\mod\,R$.
By the assumption and Lemma \ref{filt}, it suffices to show that $R/\p$ belongs to $\lap\,\G$ for every $\p\in\Spec\,R$.
Proposition \ref{laprap}(1) says that one has only to prove that $\tr (R/\p )$ belongs to $\rap\,\G$ for every $\p\in\Spec\,R$.
There is an exact sequence
$$
R^n \to R \to R/\p \to 0,
$$
where $n=\nu _R (\p )$.
Taking the $R$-dual of this sequence, we get another exact sequence
$$
0 \to (0: \p ) \to R \to R^n \to \tr (R/\p ) \to 0.
$$
Noting that $R$ is an integral domain, we see that
$$
(0:\p )=
\begin{cases}
R & \text{ if }\p =0,\\
0 & \text{ if }\p\neq 0.
\end{cases}
$$
In particular, the $R$-module $(0: \p )$ is free.
Hence the $R$-module $\tr (R/\p )$ has projective dimension at most two, in particular, it has finite G-dimension.
Corollary \ref{abg}(1) says that $\tr (R/ \p )$ belongs to $\rap\,\G$, as desired.
\qed
\end{pf}


{\it Acknowledgments.}
The author would like to thank Yuji Yoshino and an anonymous referee for helpful comments and useful suggestions.


\end{document}